
\magnification\magstephalf
\documentstyle{amsppt}

\hsize 5.72 truein
\vsize 7.9 truein
\hoffset .39 truein
\voffset .26 truein
\mathsurround 1.67pt
\parindent 20pt
\normalbaselineskip 13.8truept
\normalbaselines
\binoppenalty 10000
\relpenalty 10000
\csname nologo\endcsname 


\font\bc=cmb10
\font\tenbsy=cmbsy10

\catcode`\@=11

\def\myitem#1.{\item"(#1)."\advance\leftskip10pt\ignorespaces}

\def\qedsymbol{{\mathsurround\z@$\square$}}
\redefine\qed{\relaxnext@\ifmmode\let\next\@qed\else
  {\unskip\nobreak\hfil\penalty50\hskip2em\null\nobreak\hfil
    \qedsymbol\parfillskip\z@\finalhyphendemerits0\par}\fi\next}
\def\@qed#1$${\belowdisplayskip\z@\belowdisplayshortskip\z@
  \postdisplaypenalty\@M\relax#1
  $$\par{\lineskip\z@\baselineskip\z@\vbox to\z@{\vss\noindent\qed}}}
\outer\redefine\beginsection#1#2\par{\par\penalty-250\bigskip\vskip\parskip
  \leftline{\tenbsy x\bf#1. #2}\nobreak\smallskip\noindent}

\def\next{\let\@sptoken= }\def\next@{ }\expandafter\next\next@
\def\@futureletnext#1{\let\nextii@#1\futurelet\next\@flti}
\def\@flti{\ifx\next\@sptoken\let\next@\@fltii\else\let\next@\nextii@\fi\next@}
\expandafter\def\expandafter\@fltii\next@{\futurelet\next\@flti}

\let\zeroindent\z@
\let\savedef@\endproclaim\let\endproclaim\relax 
\define\chkproclaim@{\add@missing\endroster\add@missing\enddefinition
  \add@missing\endproclaim
  \envir@stack\endproclaim
  \edef\endit@{\leftskip\the\leftskip\rightskip\the\rightskip}}
\let\endproclaim\savedef@
\def\thing@{.\enspace\egroup\ignorespaces}
\def\thingi@(#1){ \rm(#1)\thing@}
\def\thingii@\cite#1{ \rm\@pcite{#1}\thing@}
\def\thingiii@{\ifx\next(\let\next\thingi@
  \else\ifx\next\cite\let\next\thingii@\else\let\next\thing@\fi\fi\next}
\def\thing#1#2#3{\chkproclaim@
  \ifvmode \medbreak \else \par\nobreak\smallskip \fi
  \noindent\advance\leftskip#1
  \hskip-#1#3\bgroup\bc#2\unskip\@futureletnext\thingiii@}
\let\savedef@\endproclaim\let\endproclaim\relax 
\def\endit{\endproclaim\endit@\let\endit@\undefined}
\let\endproclaim\savedef@
\def\defn#1{\thing\parindent{Definition #1}\rm}
\def\lemma#1{\thing\parindent{Lemma #1}\sl}
\def\prop#1{\thing\parindent{Proposition #1}\sl}
\def\thm#1{\thing\parindent{Theorem #1}\sl}
\def\cor#1{\thing\parindent{Corollary #1}\sl}

\def\example#1{\thing\zeroindent{Example #1}\rm}
\def\narrowthing#1{\chkproclaim@\medbreak\narrower\noindent
  \it\def\next{#1}\def\next@{}\ifx\next\next@\ignorespaces
  \else\bgroup\bc#1\unskip\let\next\narrowthing@\fi\next}
\def\narrowthing@{\@futureletnext\thingiii@}

\def\@cite#1,#2\end@{{\rm([\bf#1\rm],#2)}}
\def\cite#1{\in@,{#1}\ifin@\def\next{\@cite#1\end@}\else
  \relaxnext@{\rm[\bf#1\rm]}\fi\next}
\def\@pcite#1{\in@,{#1}\ifin@\def\next{\@cite#1\end@}\else
  \relaxnext@{\rm([\bf#1\rm])}\fi\next}

\advance\minaw@ 1.2\ex@
\atdef@[#1]{\ampersand@\let\@hook0\let\@twohead0\brack@i#1,\z@,}
\def\brack@{\z@}
\let\@@hook\brack@
\let\@@twohead\brack@
\def\brack@i#1,{\def\next{#1}\ifx\next\brack@
  \let\next\brack@ii
  \else \expandafter\ifx\csname @@#1\endcsname\brack@
    \expandafter\let\csname @#1\endcsname1\let\next\brack@i
    \else \Err@{Unrecognized option in @[}%
  \fi\fi\next}
\def\brack@ii{\futurelet\next\brack@iii}
\def\brack@iii{\ifx\next>\let\next\brack@gtr
  \else\ifx\next<\let\next\brack@less
    \else\relaxnext@\Err@{Only < or > may be used here}
  \fi\fi\next}
\def\brack@gtr>#1>#2>{\setboxz@h{$\m@th\ssize\;{#1}\;\;$}%
 \setbox@ne\hbox{$\m@th\ssize\;{#2}\;\;$}\setbox\tw@\hbox{$\m@th#2$}%
 \ifCD@\global\bigaw@\minCDaw@\else\global\bigaw@\minaw@\fi
 \ifdim\wdz@>\bigaw@\global\bigaw@\wdz@\fi
 \ifdim\wd@ne>\bigaw@\global\bigaw@\wd@ne\fi
 \ifCD@\enskip\fi
 \mathrel{\mathop{\hbox to\bigaw@{$\ifx\@hook1\lhook\mathrel{\mkern-9mu}\fi
  \setboxz@h{$\displaystyle-\m@th$}\ht\z@\z@
  \displaystyle\m@th\copy\z@\mkern-6mu\cleaders
  \hbox{$\displaystyle\mkern-2mu\box\z@\mkern-2mu$}\hfill
  \mkern-6mu\mathord\ifx\@twohead1\twoheadrightarrow\else\rightarrow\fi$}}%
 \ifdim\wd\tw@>\z@\limits^{#1}_{#2}\else\limits^{#1}\fi}%
 \ifCD@\enskip\fi\ampersand@}
\def\brack@less<#1<#2<{\setboxz@h{$\m@th\ssize\;\;{#1}\;$}%
 \setbox@ne\hbox{$\m@th\ssize\;\;{#2}\;$}\setbox\tw@\hbox{$\m@th#2$}%
 \ifCD@\global\bigaw@\minCDaw@\else\global\bigaw@\minaw@\fi
 \ifdim\wdz@>\bigaw@\global\bigaw@\wdz@\fi
 \ifdim\wd@ne>\bigaw@\global\bigaw@\wd@ne\fi
 \ifCD@\enskip\fi
 \mathrel{\mathop{\hbox to\bigaw@{$%
  \setboxz@h{$\displaystyle-\m@th$}\ht\z@\z@
  \displaystyle\m@th\mathord\ifx\@twohead1\twoheadleftarrow\else\leftarrow\fi
  \mkern-6mu\cleaders
  \hbox{$\displaystyle\mkern-2mu\copy\z@\mkern-2mu$}\hfill
  \mkern-6mu\box\z@\ifx\@hook1\mkern-9mu\rhook\fi$}}%
 \ifdim\wd\tw@>\z@\limits^{#1}_{#2}\else\limits^{#1}\fi}%
 \ifCD@\enskip\fi\ampersand@}

\define\eg{{\it e.g\.}}
\define\ie{{\it i.e\.}}
\define\today{\number\day\ \ifcase\month\or
  January\or February\or March\or April\or May\or June\or
  July\or August\or September\or October\or November\or December\fi
  \ \number\year}
\def\pr@m@s{\ifx'\next\let\nxt\pr@@@s \else\ifx^\next\let\nxt\pr@@@t
  \else\let\nxt\egroup\fi\fi \nxt}

\define\widebar#1{\mathchoice
  {\setbox0\hbox{\mathsurround\z@$\displaystyle{#1}$}\dimen@.1\wd\z@
    \ifdim\wd\z@<.4em\relax \dimen@ -.16em\advance\dimen@.5\wd\z@ \fi
    \ifdim\wd\z@>2.5em\relax \dimen@.25em\relax \fi
    \kern\dimen@ \overline{\kern-\dimen@ \box0\kern-\dimen@}\kern\dimen@}%
  {\setbox0\hbox{\mathsurround\z@$\textstyle{#1}$}\dimen@.1\wd\z@
    \ifdim\wd\z@<.4em\relax \dimen@ -.16em\advance\dimen@.5\wd\z@ \fi
    \ifdim\wd\z@>2.5em\relax \dimen@.25em\relax \fi
    \kern\dimen@ \overline{\kern-\dimen@ \box0\kern-\dimen@}\kern\dimen@}%
  {\setbox0\hbox{\mathsurround\z@$\scriptstyle{#1}$}\dimen@.1\wd\z@
    \ifdim\wd\z@<.28em\relax \dimen@ -.112em\advance\dimen@.5\wd\z@ \fi
    \ifdim\wd\z@>1.75em\relax \dimen@.175em\relax \fi
    \kern\dimen@ \overline{\kern-\dimen@ \box0\kern-\dimen@}\kern\dimen@}%
  {\setbox0\hbox{\mathsurround\z@$\scriptscriptstyle{#1}$}\dimen@.1\wd\z@
    \ifdim\wd\z@<.2em\relax \dimen@ -.08em\advance\dimen@.5\wd\z@ \fi
    \ifdim\wd\z@>1.25em\relax \dimen@.125em\relax \fi
    \kern\dimen@ \overline{\kern-\dimen@ \box0\kern-\dimen@}\kern\dimen@}%
  }

\catcode`\@\active

\let\PVstyle=d 

\font\tenscr=rsfs10 
\font\sevenscr=rsfs7 
\font\fivescr=rsfs5 
\skewchar\tenscr='177 \skewchar\sevenscr='177 \skewchar\fivescr='177
\newfam\scrfam \textfont\scrfam=\tenscr \scriptfont\scrfam=\sevenscr
\scriptscriptfont\scrfam=\fivescr
\define\scr#1{{\fam\scrfam#1}}
\let\Cal\scr

\let\0\relax 
\mathchardef\idot="202E
\define\restrictedto#1{\big|_{#1}}
\define\pr{\operatorname{pr}}
\define\Gm{\Bbb G_{\text{m}}}
\define\Ker{\operatorname{Ker}}
\define\Num{\operatorname{Num}}
\define\Pic{\operatorname{Pic}}
\define\Spec{\operatorname{Spec}}
\define\Supp{\operatorname{Supp}}

\hyphenation{pa-ram-e-trized}
\hyphenation{semi-abel-i-an}

\topmatter
\title Integral Points on Subvarieties of Semiabelian Varieties,~II\endtitle
\author Paul Vojta\snug${}^{*}$\snug\endauthor
\leftheadtext{Paul Vojta}
\affil University of California, Berkeley\endaffil
\address University of California, Department of Mathematics \#3840,
  Berkeley, CA \ 94720--3840, USA\endaddress
\date 21 May 1998\enddate
\thanks \snug${}^{*}$\snug Partially supported by NSF grants DMS-9304899
  and DMS-9532018\endthanks

\abstract
This paper proves a finiteness result for families of integral points on
a semiabelian variety minus a divisor, generalizing the corresponding
result of Faltings for abelian varieties.  Combined with the main theorem
of the first part of this paper, this gives a finiteness statement for
integral points on a closed subvariety of a semiabelian variety,
minus a divisor.

In addition, the last two sections generalize some standard results on
closed subvarieties of semiabelian varieties to the context of closed
subvarieties minus divisors.
\endabstract
\endtopmatter

\document
Recall that a semiabelian variety is a group variety $A$ such that, after
suitable base change, there exists an abelian variety $A_0$ and an exact
sequence
$$0 \to \Gm^\mu \to A @>\rho>> A_0\to 0\;.\tag\00.1$$
(In this paper a {\bc variety} is an integral separated scheme of finite type
over a field.  Since a group variety has a rational point, the base field
is algebraically closed in the function field.  In characteristic zero,
this implies that the variety is geometrically integral.)

Let $k$ be a number field, and let $S$ be a finite set of places of $k$
containing all archimedean places.  Let $R$ be the ring of integers of $k$
and let $R_S$ be the localization of $R$ away from places in $S$.
Let $X$ be a quasi-projective variety over $k$.  A {\bc model for $X$
over $R_S$} is an integral scheme, surjective and quasi-projective over
$\Spec R_S$, together with an isomorphism of the generic fiber over $k$
with $X$.  An {\bc integral point} of $\Cal X$ (or, loosely speaking,
an integral point of $X$) is an element of $\Cal X(R_S)$.

The first part \cite{V~3} of this paper proved a finiteness statement
(Theorem \00.3) for families of integral points on closed subvarieties
$X$ of a semiabelian variety $A$ over $k$.  This second and final part
proves a similar result (Theorem \00.2) for certain open subvarieties
of $A$.

Specifically, the open subvarieties under consideration are those that
can be written as the complement of a divisor.  These two finiteness
results then combine very easily to give a finiteness statement
(Theorem \00.4) for a closed subvariety minus a divisor.

The main result of this paper is the following:

\thm{\00.2}
Let $D$ be an effective divisor on $A$, let $V=A\setminus\Supp D$,
and let $\Cal V$ be a model for $V$ over $\Spec R_S$.
Then the set $\Cal V(R_S)$ of integral points on $V$ equals a
finite union $\bigcup_i\Cal B_i(R_S)$, where each $\Cal B_i$ is
a subscheme of $\Cal V$ whose generic fiber $B_i$ is a translated
semiabelian subvariety of $A$.
\endit

In \cite{F, Thm.~2}, Faltings proved a corresponding statement for integral
points on abelian varieties:  if $D$ is an ample effective divisor on an
abelian variety $A$, then (with notations as above) the set $\Cal V(R_S)$
is {\it finite.}  In the semiabelian case this is no longer true (see
Examples \01.1 and \01.3); however it is true that Theorem \00.2
implies Faltings' result.  Indeed, since an ample divisor $D$ on
an abelian variety $A$ is still ample when restricted to
a nontrivial translated abelian subvariety, the result
follows by induction on the dimension of $A$.

As with Faltings' result, the proof of Theorem \00.2 proceeds by reducing
to a statement on diophantine approximation (Theorem \03.6); in addition,
this paper relies heavily on results from \cite{V~3}.

In \cite{V~3} we proved the following result.

\thm{\00.3}
Let $X$ be a closed subvariety of a semiabelian variety $A$,
and let $\Cal X$ be a model for $X$.  Then $\Cal X(R_S)$ equals a finite
union $\bigcup_i\Cal B_i(R_S)$, where each $\Cal B_i$ is a subscheme
of $\Cal X$ whose generic fiber $B_i$ is a translated semiabelian subvariety
of $A$.
\endit

Since the conclusion of this theorem is the same as for Theorem \00.2,
the two theorems can be combined and generalized:

\thm{\00.4}
Let $X$ be a closed subvariety of a semiabelian variety $A$,
let $D$ be an effective divisor on $X$, and let $\Cal V$ be a
model for $X\setminus D$.  Then $\Cal V(R_S)$ equals a finite
union $\bigcup_i\Cal B_i(R_S)$, where each $\Cal B_i$ is a subscheme
of $\Cal V$ whose generic fiber $B_i$ is a translated semiabelian subvariety
of $A$.
\endit

\demo{Proof}
Let $\Cal X$ be a model for $X$ such that $\Cal V\subseteq\Cal X$
(this can be accomplished, if necessary, by enlarging $S$).  By Theorem \00.2,
$$\Cal V(R_S) \subseteq \Cal X(R_S) = \bigcup_i\Cal B_i(R_S)\; ;$$
then Theorem \00.3 implies that
$$\Cal V(R_S) = \bigcup_i(\Cal B_i\cap\Cal V)(R_S)
  = \bigcup_i \bigcup_j \Cal B_{ij}(R_S)\;.\qed$$
\enddemo

The situation regarding integral points on complements of sets of
codimension $\ge2$ is not as clean; in this case most of the rational
points are also integral.

I do not believe that Theorem \00.2 has been conjectured by anyone,
except that it follows from the general conjectures of \cite{V~1}
(see Theorem \05.16).

The first section of the paper gives some examples showing that a stronger
conclusion in Theorem \00.2 is impossible.  The second section begins the
proof proper by showing some results on completions of semiabelian
varieties; it is this section that contains most of what is new.
Section \03 proves the main approximation lemma via extensions of
Thue's method from \cite{F} and \cite{V~3}.  Section \04 then combines
these results to conclude the proof of Theorem \00.2.
The last two sections prove some theorems suggested by the similarities
between Theorems \00.2 and \00.3.

The author would like to thank S. Sperber for useful suggestions.
He also thanks ETH-Z\"urich and the Max-Planck-Institut in Bonn
for hospitality during brief stays.

\beginsection{\01}{Some examples}

This section gives some examples showing that one cannot hope to get
finiteness of the set of integral points in Theorem \00.2.

\example{\01.1}
Let $A=\Gm^2$ and let $D$ be the diagonal on $A$.
Completing $A$ in the obvious way to $(\Bbb P^1)^2$, it follows that
the closure of $D$ is ample.  Yet
$$A\setminus D\cong(\Gm\setminus\{1\})\times\Gm\; ,$$
so it may have infinitely many integral points (depending on the model).
\endit

Of course in this case there is a nontrivial group action on $A$.
The following definition, which will be useful throughout this paper,
formalizes this idea.

\defn{\01.2}  Let $X$ be a variety on which a semiabelian variety $A$ acts,
and let $Y$ be a subvariety of $X$.
\roster
\myitem a.  Let $B(A,Y)$ denote the identity component of the subgroup
$$\{a\in A(\widebar{\Bbb Q})\mid a+Y=Y\}$$
in $A$.
\myitem b.  If $X=A$, acting on itself by translation, then we write
$B(Y)=B(A,Y)$.  The restriction to $Y$ of the quotient map $A\to A/B(Y)$
exhibits $Y$ as a fiber bundle with fiber $B(Y)$.  This map is called the
{\bc Ueno fibration} associated to $Y$.  It is {\bc trivial} when $B(Y)$ is.
\endroster
\endit

Classically, the Ueno fibration is defined for closed subvarieties.
The image of the Ueno fibration has trivial Ueno fibration.

\example{\01.3}
Let $E$ be an elliptic curve and let $A=\Gm\times E$.
Let $f$ be a nonzero rational function on $E$ with a pole at a rational
point $P$.  Let $U$ be the largest subset of $E$ on which $f$ is defined
and nonzero, and let $D\subseteq\Gm\times U\subseteq A$ be its graph.
Then $A\setminus D$ has trivial Ueno fibration, yet it contains the
nontrivial translated subgroup $\Gm\times\{P\}$.
\endit

Both of these examples illustrate that the non-completeness of semiabelian
varieties introduces issues not present in the case of abelian varieties.

\beginsection{\02}{Completions of semiabelian varieties}

This section collects some results about completions of semiabelian
varieties with various desirable properties.

\smallskip\bgroup\narrower\it Throughout this section, all varieties
are over a field of characteristic zero, although it may well be true
that everything holds over fields of arbitrary characteristic.
\par\egroup\smallbreak

\defn{\02.1}
Let $G$ be a variety.  A {\bc completion of $G$}
is a complete variety $X$ with an open immersion $G\hookrightarrow X$.
We often identify $G$ with its image in $X$.
Given two completions $X_1$ and $X_2$ of $G$,
we say that {\bc $X_1$ dominates $X_2$} if there exists a morphism $X_1\to X_2$
compatible with the immersions $G\hookrightarrow X_1$
and $G\hookrightarrow X_2$.
A completion $X$ of a group variety $G$ is {\bc equivariant}
if the group law $G\times G\to G$ extends to a morphism $G\times X\to X$.
\endit

As noted in \cite{V~3, Sect.~2}, a semiabelian variety $A$ is isomorphic to
$$\Bbb P'(\Cal O_{A_0}\oplus\Cal M_1)\times_{A_0}\dots
  \times_{A_0}\Bbb P'(\Cal O_{A_0}\oplus\Cal M_\mu)$$
for some $\Cal M_1,\dots,\Cal M_\mu\in\Pic^0 A_0$.
Here the notation $\Bbb P'(\Cal L_1\oplus\Cal L_2)$ means
$\Bbb P(\Cal L_1\oplus\Cal L_2)$
minus the sections corresponding to the canonical projections to $\Cal L_1$
and $\Cal L_2$.  This paper needs a slightly more general situation in which $V$
is a projective variety, $\Cal M_1,\dots,\Cal M_n\in\Pic^0 V$\snug,
and $W=\Bbb P'(\Cal O_V\oplus\Cal M_1)\times_V\dots
  \times_V\Bbb P'(\Cal O_V\oplus\Cal M_\mu)$.  The group $\Gm^\mu$ still
acts on $W$, and $W$ is a fiber bundle over $V$ with fiber $\Bbb G_m^\mu$.
As usual, let $\rho\:W\to V$ denote the canonical projection.
In practice $V$ will be birational to $A_0$.

Throughout this paper, all fiber bundles will have fiber equal to the
variety underlying a group variety.  The structure group of such bundles
will always be the group of translations.

\lemma{\02.2}
Any equivariant completion $\widebar G$
of $\Gm^\mu$ determines a completion $\widebar W$ of $W$
for which $\rho$ extends to a fiber bundle $\bar\rho\:\widebar W\to V$
with fiber $\widebar G$.
\endit

\demo{Proof}
Cover $V$ by open subsets $U_i$ on which $\rho^{-1}(U_i)$ is
isomorphic to a product $\Gm^\mu\times U_i$.  We will form $\widebar W$
by glueing completions $\widebar G\times U_i$ of $\rho^{-1}(U_i)$.
This is possible since the glueing isomorphisms on $\rho^{-1}(U_i\cap U_j)$
are translations by sections of $\rho$, and $\widebar G$ is an equivariant
completion.  See also \cite{Se~2, Sect.~1.3}.\qed
\enddemo

If $V=A_0$ then $W=A$ and the resulting completion
$\widebar W=:\widebar A$ is then invariant under translation by elements of $A$.

\prop{\02.3}
Let $A$ and $\rho\:A\to A_0$ be as in (\00.1), let $A'$
be another semiabelian variety, and let $\theta\:A\to A'$ be a group
homomorphism.  Let $X$ and $X'$ be equivariant completions
of $A$ and $A'$, respectively, as in Lemma \02.2.  Let $L$ and $L'$
be ample line sheaves on $X$ and $X'$, respectively, and let $h_L$ and $h_{L'}$
be associated height functions.  Then for all $P\in A(\widebar{\Bbb Q})$,
$$h_{L'}(\theta(P)) \ll h_L(P) + O(1)\;.\tag\02.3.1$$
\endit

\demo{Proof}
Let $\Gamma$ be the closure of the graph of $\theta$
in $X\times X'$, let $M$ be an ample line sheaf on $\Gamma$, and let $h_M$
be an associated height function.  Then basic properties of heights
(functoriality and \cite{V~1, 1.2.9f}) imply that
$$h_{L'}(\theta(P)) \ll h_M((P,\theta(P))) + O(1)\;.$$
Thus, by replacing $X'$ with $\Gamma$ and $A'$ with $A$, we reduce to
the case where $A'=A$ and $X'$ dominates $X$ (since $\Gamma$ is also
an equivariant completion of $A$).
Let $\phi\:X'\to X$ be the morphism.  By Kodaira's lemma \cite{V~1, 1.2.7},
we have $m\phi^{*}L\sim L'+E$ for some $m>0$ and some effective divisor $E$
on $X'$.  Thus (\02.3.1) holds outside of the base locus $B$ of $E$.
Let $\tau\:A\to A$ denote translation by an element $a\in A$.
Then for $P\notin\tau^{-1}(B)$,
$$\split h_{L'}(\theta(P)) &\ll h_{\tau^{*}L'}(\theta(P)) + O(1) \\
  &= h_{L'}(\tau(\theta(P))) + O(1) \\
  &\ll h_L(\tau(P)) + O(1) \\
  &= h_{\tau^{*}L}(P) + O(1) \\
  &\ll h_L(P) + O(1)\;.\endsplit\tag\02.3.2$$
Here the first and last steps follow from \cite{V~1, 1.2.9d},
since $\tau^{*}L$ is algebraically equivalent to $L$ on $X$.
The proposition then follows by applying (\02.3.2) to
elements $a_1,\dots,a_r\in A$ chosen
such that the corresponding translations $\tau_1,\dots,\tau_r$ satisfy
$$\bigcap_{i=1}^r \tau_i^{-1}(B\cap A)=\emptyset\;.\qed$$
\enddemo

\thm{\02.4}
Let $A_0$ be a smooth projective variety,
let $\Cal M_1,\dots,\Cal M_\mu\in\Pic^0 A_0$, and let
$$A=\Bbb P'(\Cal O_{A_0}\oplus\Cal M_1)\times_{A_0}\dots\times_{A_0}
  \Bbb P'(\Cal O_{A_0}\oplus\Cal M_\mu)\;.$$
Let $D$ be an effective divisor on $A$.
Assume that $B(\Gm^\mu,D)=0$ (where $\Gm^\mu$ acts on fibers of $A$ over $A_0$
in the obvious manner).
Then there exists an equivariant completion $\widebar G$
of $\Gm^\mu$ with corresponding completion $\widebar A$ of $A$
(as in Lemma \02.2) and a projective birational morphism
$\pi_0\:\widetilde A_0\to A_0$, satisfying the following conditions.
\roster
\myitem 1.  Let $\widetilde A=\widebar A\times_{A_0}\widetilde A_0$ and
let $\pi\:\widetilde A\to \widebar A$ be the canonical projection onto
the first factor.  Then there exists a Cartier divisor $\widetilde D$
on $\widetilde A$ such that $(\pi_{*}\widetilde D)\cap A=D$,
where $\pi_{*}$ refers to $\widetilde D$ as a Weil divisor.
Moreover, $\widetilde D$ is ample on fibers of the map
$\tilde\rho\:\widetilde A\to A_0$.
\myitem 2.  The support of $\widetilde D$ does not contain $\pi^{-1}(T)$ for
any subset $T\subseteq\widebar A$ corresponding to a $G$\snug-orbit
of $\widebar G$.
\myitem 3.  Let $E$ be the exceptional set for $\pi_0\:\widetilde A_0\to A_0$
and let $\tilde\rho\:\widetilde A\to\widetilde A_0$ be the map corresponding
to $\bar\rho\:\widebar A\to A_0$.  Then there exists a Cartier divisor
$\widebar D$ on $\widebar A$ and a divisor $F$ on $\widetilde A$ supported only
on $\tilde\rho^{-1}(E)$ satisfying the numerical equivalence
$$\pi^{*}\widebar D \equiv \widetilde D + F\;.\tag\02.4.1$$
\myitem 4.  The set $\rho^{-1}(\pi_0(E))$ is contained in $\Supp D$.
\myitem 5.  Finally, $\widebar G$ has only finitely many $G$\snug-orbits.
\endroster
Moreover, if $A_0$ is abelian (and $A$ is semiabelian) and if $D$ has
trivial Ueno fibration, then $\widebar D$ is ample.
$$\CD \widetilde A @>\tilde\rho>> \widetilde A_0 \\
  @V\pi VV @V\pi_0 VV \\
  \widebar A @>\bar\rho>> A_0 \\
  \cup& \nearrow\scriptstyle{\rho} \\
  A\endCD$$
\endit

\demo{Proof}
We first consider the case $A=\Gm^\mu$.  In this case $\pi$
is an isomorphism, (4) is vacuous, and conditions (1) and (3) are
equivalent to the closure $\widebar D$ of $D$ in $\widebar A$ being Cartier
and ample.

Let $f\in k[x_1,x_1^{-1},\dots,x_\mu,x_\mu^{-1}]$ be a Laurent polynomial
such that $D=(f)$.  Let $M=\Bbb Z^\mu$, and for $m=(m_1,\dots,m_\mu)\in M$
let $x^m$ denote $x_1^{m_1}\dotsm x_\mu^{m_\mu}$.  Write $f=\sum c_mx^m$
and let $\Delta=\Delta_f$ be the convex hull
of $\{m\in M\mid c_m\ne0\}$ in $M_{\Bbb R}:=M\otimes_{\Bbb Z}\Bbb R$.
This is a polyhedron called the {\bc Newton polyhedron.}
Its vertices are lattice points.
Let $m^{(0)},\dots,m^{(\ell)}$ be the lattice points in $\Delta$.
For vertices $m$ of $\Delta$, let $\sigma_m\subseteq M_{\Bbb R}$ be the cone
$$\Biggl\{\sum_{i=0}^\ell \lambda_i(m^{(i)}-m)
  \Bigm| \lambda_i\ge0\ \forall\ i\Biggr\}\;.$$
After replacing $\Delta$ with a positive integral multiple $n\Delta$,
we may assume that for all vertices $m$ of $\Delta$,
the set $\{m^{(i)}-m\mid i=0,\dots,\ell\}$ generates
the monoid $\sigma_m\cap M$ (see Gordan's lemma, \cite{TE, p\. 7}).
This multiple $n\Delta$ corresponds to $f^n$, which corresponds to $nD$.

Following \cite{O, Sect.~1},
we have a morphism $\phi\:\Gm^\mu\to\Bbb P^\ell$ defined by
$$(t_1,\dots,t_\mu)\mapsto[t^{m^{(0)}}:\dots:t^{m^{(\ell)}}]\;.$$
Since $D$ has trivial Ueno fibration, the Newton polyhedron does not
lie in any hyperplane of $M_{\Bbb R}$; therefore this map is actually an
embedding.  Let $\widebar G$ be the closure of the image of $\phi$.
Then $\widebar G$ is a {\bc toric variety} (for definitions see \cite{TE}
or \cite{D}).
In particular it is an equivariant completion of $\Gm^\mu$
with only finitely many orbits.  These orbits are in
incidence-preserving one-to-one correspondence with the faces of $\Delta$.
Finally, by \cite{TE, p\. 6 Thm.~1}, $\widebar G$ is normal.

\lemma{\02.4.2}
Let $D$ be a divisor on $\Gm^\mu$ with trivial Ueno fibration.
Let $\widebar G$ be the toric variety corresponding to a suitably large
multiple of the Newton polyhedron $\Delta$ of the defining equation $f$ of $D$.
Then the closure $\widebar D$ of $D$ in $\widebar G$ is Cartier and ample,
and its support does not contain any $\Gm^\mu$\snug-orbit of $\widebar G$.
\endit

\demo{Proof}
The polynomial $f$ defines a global section of $\Cal O(1)$
on $\Bbb P^\ell$ by $s=c_0x_0+\dots+c_\ell x_\ell$\snug,
where $f=\sum c_i x^{m^{(i)}}$.
Given a face $\delta$ of $\Delta$, the closure of the associated orbit
in $\widebar G$ is determined by the vanishing of all $x_i$
for which $m^{(i)}\notin\bar\delta$.  By definition of $\delta$,
there is an index $i$ for which $c_i\ne0$ and $m_i\in\bar\delta$;
therefore $s$ does not vanish identically on the orbit associated to $\delta$.
Thus $(s)=\widebar D$, which is therefore ample.\qed
\enddemo

The lemma implies conditions (1)--(4), thus proving the case $A=\Gm^\mu$.
Before proving the general case, some more lemmas are needed.

\lemma{\02.4.3}
Let $\widebar G$ be a completion of $\Gm^\mu$ as described above,
and let $V$, $W$, and $\widebar W$ be as in Lemma \02.2.
Recall that $\Num X$ denotes the numerical equivalence class group of
a variety $X$.  Let $i_1\:\widebar G\to\widebar W$ be a closed fiber of
$\bar\rho\:\widebar W\to V$ (equivariant under the action of $\Gm^\mu$),
and let $i_2\:V\to\widebar W$ be a section of $\bar\rho$ associated to
a zero-dimensional orbit of $\widebar G$.  Then
\roster
\myitem a.  the map
$$(i_1^{*},i_2^{*})\:\Num\widebar W \to \Num\widebar G\times\Num V$$
is an isomorphism, and is independent of the choice of $i_1$ and $i_2$;
\myitem b.  every closed integral curve on $\widebar W$ is numerically
equivalent to an effective sum of curves in the images of $i_1$ and $i_2$; and
\myitem c.  a divisor $D$ on $\widebar W$ is ample if and only if the
divisors $i_1^{*}D\in\Pic\widebar G$ and $i_2^{*}D\in\Pic V$ are ample.
\endroster
\endit

\demo{Proof}  First consider part (b).  We start by claiming that every
closed integral curve on $\widebar W$ is numerically equivalent to an
effective sum of curves in fibers of $\bar\rho$ and curves in sections
of $\bar\rho$ associated to zero-dimensional orbits of $\widebar G$.
Let $C$ be a closed integral curve in $\widebar W$ and let
$\chi\:\Gm\to\Gm^\mu$ be a nontrivial one-parameter subgroup.
Translations of $C$ by $\chi(a)$ for $a\in\Gm$ define a surface
$Y\subseteq\widebar W\times\Gm$ with $Y\cap(\widebar W\times\{1\})=C$.
Let $\widebar Y$ be the closure of $Y$ in $\widebar W\times\Bbb P^1$,
and let $Y_0=\widebar Y\cap(\widebar W\times\{0\})$.  Then $Y_0$ is a sum of
curves in $\widebar W$ which is numerically equivalent to $C$ and is
invariant under translations by $\chi$.  Thus each irreducible component
either lies in $\widebar W\setminus W$ or lies in a fiber
of $\bar\rho\:\widebar W\to V$.  The claim then follows by induction on
dimension.

Given any two fibers of $\bar\rho$, an effective sum of curves in one fiber
is numerically equivalent to an effective sum of curves in the other.
Hence, in the above claim, the curves in fibers of $\bar\rho$ can be
assumed to lie in the image of $i_1$.

To prove (b), it remains to show that the curves in sections of $\bar\rho$,
as above, can be taken to lie in the image of $i_2$.  To show this, it suffices
to show that for any two sections $\sigma_1$ and $\sigma_2$ of $\bar\rho$
as above and any closed integral curve $C\subseteq V$, $\sigma_1(C)$ is
algebraically equivalent to $\sigma_2(C)$.
Let $\Gamma$ be the graph whose vertices are zero-dimensional orbits of
$\widebar G$ and whose edges are one-dimensional orbits.  Since $\Gamma$
is connected, it suffices to consider sections $\sigma_1$, $\sigma_2$ contained
in the subset of $\widebar W$ corresponding to the closure of a one-dimensional
orbit in $\widebar G$.  In this case it is easy to show explicitly that
$\sigma_1(C)$ and $\sigma_2(C)$ are algebraically equivalent, since this
subset of $\widebar W$ is isomorphic to $\Bbb P(\Cal O_V\oplus\Cal M)$
for some $\Cal M\in\Pic^0 V$.  Thus (b) holds.  Moreover,
$i_2^{*}\:\Num\widebar W\to\Num V$ is independent of the choice of $i_2$.

Next consider the map
$$i_1^{*}\:\Pic\widebar W\to\Pic\widebar G\;.$$
Since $\widebar G$ has trivial Albanese, this map is independent of the
fiber chosen.  Since divisors on $\Gm^\mu$ are all principal, it follows
that $\Pic\widebar G$ is generated by the closures of the orbits of
codimension one; hence $i_1^{*}\:\Pic\widebar W\to\Pic\widebar G$ is surjective.
By the Seesaw theorem \cite{Mi, Thm.~5.1}, the kernel is $\Pic V$.
Thus there is an exact sequence
$$0 \to \Pic V \to \Pic\widebar W \to \Pic\widebar G \to 0\;.$$
The map $i_2^{*}\:\Pic\widebar W\to\Pic V$ splits this sequence.
(This splitting depends on the choice of $i_2$.)

By part (b), an element in $\Pic\widebar W$ is numerically equivalent
to zero if and only if its components in $\Pic V$ and $\Pic\widebar G$
are both numerically equivalent to zero.  Hence there is an exact sequence
$$0 \to \Num V \to \Num\widebar W \to \Num\widebar G \to 0\; ,$$
which is again split.  In this case, though, the splitting is independent
of the choice of zero-dimensional orbit, by the argument using $\Gamma$.

Since $i_1$ and $i_2$ are closed immersions, it follows that if $D$ is
an ample divisor on $\widebar W$ then its components in $\Num\widebar G$
and $\Num V$ must also be ample.  The converse follows from part (b)
and from Kleiman's criterion for ampleness.\qed
\enddemo

\lemma{\02.4.4}
Let $A$ be a semiabelian variety with $\mu=1$,
let $\widebar A$ be the completion of $A$ associated to
the (unique) completion $\Gm\subseteq\Bbb P^1$, and
let $\bar\rho\:\widebar A\to A_0$ be the extension of $\rho\:A\to A_0$
as in (\00.1).  Let $\sigma\:A_0\to\widebar A$ and $\sigma'\:A_0\to\widebar A$
be the sections of $\bar\rho$ corresponding to $0\in\Bbb P^1$
and $\infty\in\Bbb P^1$, respectively.
Let $D$ be a closed subset of $\widebar A$ of pure codimension one.
Assume that $D$ meets the generic fiber of $\rho$, but that it does not
contain the image of $\sigma$ or $\sigma'$.
Then (in the notation of Definition \01.2) there exists
an abelian subvariety $C$ of $A$ such that
$$\rho(C) = B(\sigma^{-1}(D))\;.$$
\endit

\demo{Proof}
By replacing $A$ with a translate of $\rho^{-1}(B(\sigma^{-1}(D)))$,
we may assume that $B(\sigma^{-1}(D))=A_0$.
We may also assume that $D$ is irreducible (discard all but one suitable
irreducible component).  Then $D$ does not meet the image of $\sigma$;
since the images of $\sigma$ and $\sigma'$ are algebraically equivalent,
it does not meet the image of $\sigma'$, either.  Thus $D\subseteq A$.
Recall that $\widebar A\cong\Bbb P(\Cal O_{A_0}\oplus\Cal M)$ for some
line sheaf $\Cal M\in\Pic^0 A_0$.  Let $Y$ be some normal projective variety
admitting a birational morphism $i\:Y\to D$.  Let $E=\widebar A\times_{A_0}Y$
and $\Cal N=i^{*}\rho^{*}\Cal M$, so that $E=\Bbb P(\Cal O_Y\oplus\Cal N)$.
This contains a divisor $D'$ (equal to the image of the map
$i\times_{A_0}\operatorname{Id}_Y$) which has degree $1$ over $Y$ and
does not meet the divisors on $E$ corresponding to the projections of
$\Cal O_Y\oplus\Cal N$ onto either of its direct summands.
Therefore the surjection $\Cal O_Y\oplus\Cal N\twoheadrightarrow\Cal L$
corresponding to $D'$ is an isomorphism on each direct summand;
thus $\Cal N\cong\Cal O_Y$.  Therefore $\Cal M$ lies in the kernel of
$\Pic^0 A_0\to\Pic^0 Y$.

But since the Albanese of $Y$ maps onto $A_0$,
the above map on $\Pic^0$ must be finite; hence $\Cal M$ is torsion
in $\Pic^0 A_0$, of order dividing the degree of $\rho\restrictedto D$.
It then follows that $D$ is a translated subgroup of $A$.
Let $C$ be that subgroup.\qed
\enddemo

\lemma{\02.4.5}
Let $\widebar G$ be the completion of $\Gm^\mu$
associated to some polyhedron $\Delta$ and let $\widebar A$ be the corresponding
completion of $A$ (Lemma \02.2).  Let $\sigma\:A_0\to\widebar A$ be the
section of $\bar\rho\:\widebar A\to A_0$ corresponding to some zero-dimensional
orbit of $\widebar G$.  Let $D$ be a closed subset of $\widebar A$
of pure codimension $1$.  Assume that $D$ meets the generic fiber
of $\bar\rho$ on all subsets of $\widebar A$ corresponding to closures
of positive dimensional orbits of $\widebar G$,
but that $D$ does not contain any subset of $\widebar A$
corresponding to a zero dimensional orbit.  Then
$$\rho(B(D\cap A)) = B(\sigma^{-1}(D))\;.$$
\endit

\demo{Proof}  Let $B=B(\sigma^{-1}(D))$.  The section $\sigma$ corresponds
to a vertex $m$ of $\Delta$.  Let $\{Y_i\mid i\in I\}$ be the set of
$A$\snug-orbits of $\widebar A$ corresponding to edges (\ie, one-dimensional
faces) of $\Delta$ incident to $m$.  The closures $\widebar Y_i$
contain the image of $\sigma$.  Applying Lemma \02.4.4 to $D\cap\widebar Y_i$
on each $Y_i$ gives an {\it abelian} subvariety mapping onto $B$.
Since the edges of $\Delta$ incident to $m$ do not
lie in any hyperplane of $M$, this implies that the line sheaves
$\Cal M_1,\dots,\Cal M_\mu$ used in defining $A$ \cite{V~3, Sect.~2}
are torsion when restricted to $B$.
Thus there exists an abelian subvariety $C\subseteq A$
for which $\rho(C)=B$.

We now show that $B(D\cap A)$ contains $C$.  After replacing $A_0$ with an
isogenous abelian variety, we may assume that $C$ has degree $1$ over $B$,
and that $\Cal M_1,\dots,\Cal M_\mu$ are trivial on $B$.  Then
$\Cal M_1,\dots,\Cal M_\mu$ all descend to line sheaves on $A_0/B$;
this then defines a fiber bundle $\widebar A\to \widebar A'$ with fiber $B$.
Since this morphism has at least one fiber which does not meet $D$,
it follows that $B(D\cap A)\supseteq C$.
Thus $\rho(B(D\cap A))\supseteq B$.  The opposite inclusion is trivial.\qed
\enddemo

We now finish the proof of Theorem \02.4.  The divisor $D$ on the generic
fiber of $\bar\rho$ defines a hyperplane in $\Bbb P^\ell$ defined over $K(A_0)$,
and therefore a rational section $\phi\:A_0\to\Bbb P(\Cal E\spcheck)$
for some appropriate vector sheaf $\Cal E$ of rank $\ell+1$ on $A_0$.
Let $\widetilde A_0$ be the closure of the graph of this rational map,
and let $\pi_0\:\widetilde A_0\to A_0$ be the canonical projection.
Then $\phi$ extends to a morphism
$\tilde\phi\:\widetilde A_0\to\Bbb P(\Cal E\spcheck)$.
Let $\widebar G$ be the completion of $\Gm^\mu$ associated to $D$
on the generic fiber of $\rho$, as above.  Let $\widebar A$ be the
associated completion of $A$; this
can be naturally identified with a subscheme of $\Bbb P(\Cal E)$,
and therefore $\tilde\phi$ defines a Cartier divisor $\widetilde D_0$
on $\widetilde A:=\widebar A\times_{A_0}\widetilde A_0$.  Moreover,
by construction no fiber of $\tilde\rho$ is contained in $\widetilde D_0$,
so since $\widetilde D_0$ is of pure codimension one, it follows that all
components of $\widetilde D_0$ map onto $\widetilde A_0$.  Thus $\widetilde D_0$
is the strict transform of the horizontal (over $A_0$) part of $D$.
Add to this the pull-back of the vertical part of $D$ to get $\widetilde D$.
This gives (1), (2), and (5).

Let $\sigma\:A_0\to\widebar A$ be a section of $\bar\rho$ as in Lemma \02.4.5
and let $\tilde\sigma\:\widetilde A_0\to\widetilde A$ be the corresponding
section of $\tilde\rho$.
Then $\widetilde D$ corresponds via Lemma \02.4.3 to the divisor
$\tilde\sigma^{*}\widetilde D$ on $\widetilde A_0$ and some ample
divisor $D_1$ on $\widebar G$.  Form $\pi_{0{*}}\tilde\sigma^{*}\widetilde D$
as a Weil divisor; it is also Cartier since $A_0$ is smooth.
Let $\widebar D$ be a divisor on $\widebar A$ corresponding to
$\pi_{0{*}}\tilde\sigma^{*}\widetilde D$ and $D_1$.  By construction,
$\pi_0^{*}(\pi_{0*}\tilde\sigma^{*}\widetilde D)-\tilde\sigma^{*}\widetilde D$
is a divisor supported only on $E$.  Pulling back to $\widetilde A$
then gives (3).

Moreover, the support of $\pi_{0{*}}\tilde\sigma^{*}\widetilde D$ is
just $\sigma^{-1}(\widebar{\Supp D})$, so
if $A_0$ is an abelian variety and if $D$ has trivial Ueno fibration,
then Lemma \02.4.5 implies that
$\pi_{0{*}}\tilde\sigma^{*}\widetilde D$ has trivial Ueno fibration.
Thus by \cite{Mu, \S6, Application 1, p\. 60},
$\pi_{0*}\tilde\sigma^{*}\widetilde D$ is ample and therefore Lemma \02.4.3
implies that $\widebar D$ is ample.

We now show (4).  Zariski's Main Theorem \cite{Ha, III 11.4} and its proof
imply that there is a Zariski-open subset $U$ of $A_0$ over which $\pi_0$
is an isomorphism, and the fibers over all $P\notin U$ are positive dimensional.
For those $P$, it follows that $\widetilde D$ maps onto $\rho^{-1}(P)$.
This gives (4).\qed
\enddemo

\prop{\02.5}
Let $\widebar A$ be the equivariant completion
of a semiabelian variety $A$ associated to some Newton polyhedron $\Delta$,
and let $T$ be an orbit of $\widebar A$.  Then there exists an open subset $U$
of $\widebar A$ containing $A$ and $T$, and an equivariant projection
$p\:U\to T$ whose restriction to $T$ is the identity.
\endit

\demo{Proof}
We continue using the notation of the proof of Theorem \02.4.

It suffices to prove the proposition in the case where $A=\Gm^\mu$.
Let $\delta$ be the face of $\Delta$ corresponding to $T$.
Let $\sigma_\delta$ be the cone in $M_{\Bbb R}$ generated by elements $m'-m$
with $m'\in\Delta$ and $m\in\delta$.  Then the set
$$U:=\Spec k[\sigma_\delta\cap M]$$
is an open affine subset of $\widebar A$.  (Here $k[\sigma_\delta\cap M]$ is a
monoid ring.)  Let $\tau_\delta$ be the largest subgroup of $\sigma_\delta$.
Then $\Spec k[\tau_\delta\cap M]=T$, and the injection $T\hookrightarrow U$
corresponds to
the surjection $k[\sigma_\delta\cap M]\twoheadrightarrow k[\tau_\delta\cap M]$
defined by $x^m\mapsto0$ for all $m\notin\tau_\sigma$.

We then let $p$ be the morphism corresponding to the canonical
injection\break
$k[\tau_\delta\cap M]\hookrightarrow k[\sigma_\delta\cap M]$.\qed
\enddemo

\beginsection{\03}{The main approximation result}

This section proves the main approximation result used in the proof
of Theorem \00.2.  First we recall a standard definition and a definition
from \cite{V~3, 7.1}.

\defn{\03.1}
A line sheaf or Cartier divisor $L$ on a complete variety $X$
is {\bc nef} (``numerically effective'') if $(L\idot C)\ge0$ for all
integral curves $C$ on $X$.
\endit

\defn{\03.2}
Let $X/k$ be a variety.  A {\bc generalized Weil
function} is a function $g\:\coprod_v U(\bar k_v)\to\Bbb R$, where $U$ is
a nonempty Zariski-open subset of $X$, such that there exists a
proper birational morphism $\Phi\:X^{*}\to X$ such that $g\circ\Phi$ extends
to a Weil function for some divisor $D^{*}$ on $X^{*}$.  It is called
{\bc effective} if $D^{*}$ is an effective divisor.  The {\bc support of
$g$}, written $\Supp g$, is defined to be the set $\Phi(\Supp D^{*})$.
N.B.:  This is {\it not\/} the set where $g\ne0$ (the analysts' definition
of support).
\endit

\defn{\03.3}
Let $Y$ be a proper closed subset of a variety $X/k$.
Then a {\bc logarithmic distance function for $Y$} is an effective
generalized Weil function $g$ on $X$ such that, if $\Phi\:X^{*}\to X$
is as in Definition \03.2, then the divisor $D^{*}$ (as above)
has $\Supp D^{*}=\Phi^{-1}(Y)$.
\endit

Note that this is not really minus the logarithm of the distance to $Y$,
especially near singularities, but we do have the following easy fact.

\lemma{\03.4}
Let $Y$ be a proper closed subset of a complete variety $X/k$,
let $\phi\:X'\to X$ be a proper birational morphism, and let $Y'=\phi^{-1}(Y)$.
Let $g$ and $g'$ be logarithmic distance functions for $Y$ and $Y'$,
respectively.  Then
$$g + O(1) \gg\ll g' + O(1)\;.$$
Moreover, $O(1)$ refers to $M_k$\snug-constants, as in \cite{L, Ch.~10, \S\,1}.
\endit

\demo{Proof}  We may assume that $X^{*}$ is the same for $g$ and $g'$;
let $D$ and $D'$ be the divisors on $X^{*}$ associated to $g$ and $g'$,
respectively.  Then since $\Supp D=\Supp D'$, it follows that $nD-D'$
and $n'D'-D$ are effective for sufficiently large $n$ and $n'$.
This implies the lemma.\qed
\enddemo

For future reference, we also note that if $X$ is an equivariant
completion of $A$, if $\lambda_\infty$ is a logarithmic distance
function for $X\setminus A$, if $X\setminus A=\bigcup\widebar T_i$,
where $T_i$ are finitely many subsets of $X$ (\eg, orbits), and if $\lambda_i$
are logarithmic distance functions for $\widebar T_i$ on $X$ for all $i$,
then $\sum\lambda_i$ is a logarithmic distance function for $X\setminus A$
and therefore
$$\lambda_\infty + O(1) \gg\ll \sum\lambda_i + O(1)\;.\tag\03.5$$
Here again $O(1)$ refers to $M_k$\snug-constants.

\thm{\03.6}
Let $X$ be an equivariant completion of a
semiabelian variety $A$.  Let $h_L$ be a height function on $X$ with respect
to an ample line sheaf $L$.
Let $\lambda_w$ be a generalized Weil function on $X$ at a place $w\in S$,
and let $\lambda_{\infty,w}$ be a logarithmic distance function
for $X\setminus A$ on $X$.
Then there do not exist a real number $\kappa>0$ and
a subset $\Cal S\subseteq A(R_S)$ satisfying the conditions (1)
$$\lambda_w(P) \ge \kappa h_L(P)\tag\03.6.1$$
for all $P\in\Cal S$; and (2) for all $\eta>0$ the inequality
$$\lambda_{\infty,w}(P) \le \eta h_L(P)\tag\03.6.2$$
holds for infinitely many $P\in\Cal S$.
\endit

\demo{Proof}
First we claim that the theorem is independent
of the choice of completion of $A$.  Indeed, suppose $X'$ is
another completion, with ample line sheaf $L'$ and
height function $h_{L'}$.
Without loss of generality we may assume that $X'$ dominates $X$
via $\phi\:X'\to X$.  By Proposition \02.3,
$$h_L(P) \gg\ll h_{L'}(P)\tag\03.6.3$$
for almost all $P\in\Cal S$.  Moreover $\lambda_w$ is also a generalized
Weil function on $X'$.  Thus (\03.6.1) holds for $X$ if and only if
it holds for $X'$ (with a different $\kappa$).  Also, Lemma \03.4 and (\03.6.3)
imply that (\03.6.2) holds for $X$ if and only if it holds for $X'$
(with a different $\eta$).

Thus we may assume that $X$ is the equivariant completion
associated to the injection $\Gm^\mu\hookrightarrow(\Bbb P^1)^\mu$.
Following \cite{V~3}, we will denote $X$ by $\widebar A$ from now on.
We may assume that $\lambda_w$ is effective.  Every generalized Weil
function is dominated by a Weil function (of a divisor on $\widebar A$),
so we may also assume that $\lambda_w=\lambda_{D,w}$
is a Weil function for an effective divisor $D$ on $\widebar A$.

The basic idea of the remainder of the proof is to incorporate
the extra machinery of \cite{V~3} into the argument of \cite{F, Sect.~6}.
As in \cite{V~3}, we let $L_0$ be an ample symmetric divisor on $A_0$,
and let $L_1=\widebar A\setminus A$ (taking all components with multiplicity
one).  Then, by basic properties of height functions, we may assume
that
$$L = \rho^{*}L_0+L_1\;.$$
Unlike \cite{F}, it is not necessary here to assume that $D$ is ample;
instead, let $\ell$ be a positive integer such that $\ell L-D$ is ample.

Let $\delta>0$ be a rational number, and choose
a positive rational $\epsilon<1$ and a positive integer $n$ satisfying
$$n\epsilon < \frac{\kappa\delta}{[k:\Bbb Q]}\tag\03.6.4$$
and
$$\frac{2\delta^n}{n!} < \frac{\epsilon^{\dim A}}{(5^{\dim A}\ell\dim A)^n}\;.
  \tag\03.6.5$$
As in \cite{V~3, Sect.~3}, let $\bold s=(s_1,\dots,s_n)$ be a tuple of
positive integers.  For integers $i$ and $j$ in $\{1,\dots,n\}$ let
$s_i\cdot\pr_i-s_j\cdot\pr_j$ denote the morphism from $A^n$ to $A$
defined using the group law.  Also as in \cite{V~3}, given any product,
$\pr_i$ denotes the projection morphism from that product to its $i^{\text{th}}$
factor.

For closed subvarieties $X_1,\dots,X_n$ of $A$,
let $\widebar X_1,\dots,\widebar X_n$ denote their respective closures
in $\widebar A$.
Let $\psi_{\bold s}\:\prod\widebar X_i\dashrightarrow\widebar A^{n-1}$
be the rational map whose components are the restrictions
of $s_i\cdot\pr_i-s_{i+1}\cdot\pr_{i+1}$ as $i$ varies from $1$ to $n-1$.
Let $W_{\bold s}$ be the closure of the graph of this rational map,
and let $\pi_{\bold s}\:W_{\bold s}\to\prod\widebar X_i$ be the natural
projection.

For $n$\snug-tuples $\bold s$ of positive integers and for rational $\epsilon$
we define
$$L_{\epsilon,\bold s}
  = \sum_{i=2}^n (s_{i-1}\cdot\pr_{i-1}-s_i\cdot\pr_i)^{*}\rho^{*}L_0
  + \sum_{i=2}^n (s_{i-1}^2\cdot\pr_{i-1}-s_i^2\cdot\pr_i)^{*}L_1
  + \epsilon\sum_{i=1}^n s_i^2\pr_i^{*}L\tag\03.6.6$$
as a $\Bbb Q$\snug-divisor class on $W_{\bold s}$ and
$$\split M_{\epsilon,\bold s}
  &= \sum_{i=2}^n(s_{i-1}\cdot\pr_{i-1}-s_i\cdot\pr_i)^{*}\rho^{*}L_0 \\
  &\qquad + s_1^2\pr_1^{*}L_1 + 2\sum_{i=2}^{n-1} s_i^2\pr_i^{*}L_1
    + s_n^2\pr_n^{*}L_1
  + \epsilon\sum_{i=1}^n s_i^2\pr_i^{*}L\;.\endsplit\tag\03.6.7$$
as a $\Bbb Q$\snug-divisor class on $\prod\widebar X_i$.
Note that these differ from their counterparts in \cite{V~3}:
the first two sums are taken over a smaller set of pairs of indices,
in line with \cite{F}.
Adding all pairs of indices is possible, but more complicated.

As in \cite{V~3}, these definitions extend by homogeneity to tuples
$\bold s$ of positive rational numbers:  let $a$ be the lowest common
denominator and let $W_{\bold s}=W_{a\bold s}$,
$\pi_{\bold s}=\pi_{a\bold s}$,
$L_{\epsilon,\bold s}=a^{-2}L_{\epsilon,a\bold s}$,
and $M_{\epsilon,\bold s}=a^{-2}M_{\epsilon,a\bold s}$.

The natural injection of \cite{V~3, 3.6} carries over to this case:
for any positive integer $d$ canceling all denominators of $\bold s$ and
$\epsilon$, we have
$$\Cal O(dL_{\epsilon,\bold s})\hookrightarrow\Cal O(dM_{\epsilon,\bold s})\;.
  \tag\03.6.8$$

With these choices, the overall strategy is the same as in
\cite{V~3, Sect.~4}, except that there is no set $Z$.
As in {\it op\. cit.,} we let $h(X_i)$ denote the height of the closed
subvariety $\widebar X_i$ of $\widebar A$, taken relative to the ample
line sheaf $L$.  We omit the subscript $L$ since heights of subvarieties
will not be taken relative to any other line sheaf (for points, however,
we will retain the subscript since heights of points will be taken relative
to other line sheaves).

The strategy is to choose $P_1,\dots,P_n\in\Cal S$ satisfying the following
conditions:
\roster
\setbox0\hbox{\bc\03.6.9.1}\rosteritemwd=\wd0
\myitem{{\bc\03.6.9.1}}.  $h_L(P_1)\ge c_1$.
\myitem{{\bc\03.6.9.2}}.  $h_L(P_{i+1})/h_L(P_i)\ge c_2\ge1$
for all $i=1,\dots,n-1$.
\myitem{{\bc\03.6.9.3}}.  $P_1,\dots,P_n$ all point in roughly the same
direction in $A(R_S)\otimes_{\Bbb Z}\Bbb R$, up to a factor $1-\epsilon_1$:
see \cite{V~3, 13.2 and 13.3}.
\endroster
The main part of the proof involves closed subvarieties $X_1,\dots,X_n$ of $A$.
We start with $X_1=\dots=X_n=A$ and successively find collections
with $\sum\dim X_i$ strictly smaller.  At each stage, $X_1,\dots,X_n$ are
assumed to satisfy the following conditions:
\roster
\setbox0\hbox{\bc\03.6.10.1}\rosteritemwd=\wd0
\myitem{{\bc\03.6.10.1}}.  Each $X_i$ contains $P_i$.
\myitem{{\bc\03.6.10.2}}.  Each $X_i$ is geometrically irreducible and defined
over $k$.
\myitem{{\bc\03.6.10.3}}.  The degrees $\deg\widebar X_i$ satisfy
$\deg\widebar X_i\le c_3$.
\myitem{{\bc\03.6.10.4}}.  The heights $h(X_i)$ will be bounded by the formula
$$\sum_{i=1}^n \frac{h(X_i)}{h_L(P_i)} \le c_4\sum_{i=1}^n \frac1{h_L(P_i)}\;.$$
\endroster
Eventually, this inductive process reaches the point where some $X_j$ is
zero dimensional; \ie, $X_j=P_j$.  When that happens, by (\03.6.10.4),
$$1 = \frac{h(X_j)}{h_L(P_j)} \le \sum_{i=1}^n\frac{h(X_i)}{h_L(P_i)}
  \le c_4\sum_{i=1}^n \frac1{h_L(P_i)} \le \frac{c_4n}{h_L(P_1)}\;.$$
This contradicts (\03.6.9.1) if $c_1>c_4n$.

Here and throughout the proof, constants $c$ and $c_i$ depend only
on $\Cal A$, $D$, $n$, $k$, $S$, $L$, and sometimes the
tuple $(\dim X_1,\dots,\dim X_n)$, but not on $P_i$, $X_i$, or $\bold s$.
The value of $c$ may change from line to line.
A precise logical statement of the inductive step appears at the
end of \cite{V~3, Sect.~4}.

For $i=1,\dots,n$ let $s_i$ be rational numbers close
to $1\big/\sqrt{h_L(P_i)}$ and let $d$ be a large sufficiently divisible
integer.  Let $\Cal V$ be the model for $\prod\widebar X_i$ constructed in
\cite{V~3, Sect.~10}; briefly, this is a model just large enough to extend
the definition of $M_{\epsilon,\bold s}$.
Let $\Gamma_\delta\bigl(\prod\widebar X_i,dM_{\epsilon,\bold s}\bigr)$
denote the submodule
of $\Gamma\bigl(\prod\widebar X_i,dM_{\epsilon,\bold s}\bigr)$
consisting of sections having index $\ge\delta$ along $D\times\dots\times D$,
relative to multiplicities $ds_1^2,\dots,ds_n^2$.
Also let $h^0_\delta\bigl(\prod\widebar X_i,dM_{\epsilon,\bold s}\bigr)$
denote the rank of this module.  We start by obtaining two estimates for
some ranks.  Except for notation, this follows \cite{F}.

\lemma{\03.6.11}
If $d$ is sufficiently divisible, then
$$h^0(W_{\bold s},dL_{\epsilon,\bold s})
  \ge \epsilon^{\dim X_1}\prod_i (ds_i^2)^{\dim X_i}
    \cdot \frac{\deg\widebar X_i}{(\dim X_i)!}
    - o\bigl(d^{\sum\dim X_i}\bigr)\tag\03.6.11.1$$
where the implicit constant in $o(\;)$ is independent of $d$.
\endit

\demo{Proof}
By homogeneity, we may assume that $s_1,\dots,s_n$ are all
integers.  We have
$$\split
  &\Biggl((s_1^2\pr_1^{*}L)^{\dim X_1} \\
  &\qquad\quad \cdot \prod_{i=2}^n
    \Bigl((s_{i-1}\cdot\pr_{i-1}-s_i\cdot\pr_i)^{*}\rho^{*}L_0
      + (s_{i-1}^2\cdot\pr_{i-1}-s_i^2\cdot\pr_i)^{*}L_1\Bigr)^{\dim X_i}
    \Biggr) \\
  &\qquad = s_1^{2\dim X_1}(\deg\widebar X_1)
    \Biggl((s_2^2\pr_2^{*}L)^{\dim X_2} \\
  &\qquad\quad \cdot \prod_{i=3}^n
    \Bigl((s_{i-1}\cdot\pr_{i-1}-s_i\cdot\pr_i)^{*}\rho^{*}L_0
      + (s_{i-1}^2\cdot\pr_{i-1}-s_i^2\cdot\pr_i)^{*}L_1\Bigr)^{\dim X_i}
    \Biggr) \\
  &\qquad = \dots = \prod_{i=1}^n s_i^{2\dim X_i}(\deg\widebar X_i)\;.
  \endsplit\tag\03.6.11.2$$
Since each of the terms
$(s_{i-1}\cdot\pr_{i-1}-s_i\cdot\pr_i)^{*}\rho^{*}L_0$,
$(s_{i-1}^2\cdot\pr_{i-1}-s_i^2\cdot\pr_i)^{*}L_1$, and $s_i^2\pr_i^{*}L$
in the definition (\03.6.6) of $L_{\epsilon,\bold s}$ is nef, (\03.6.11.2)
gives
$$\left(L_{\epsilon,\bold s}^{\sum\dim X_i}\right)
 \ge \epsilon^{\dim X_1}\binom{\sum\dim X_i}{\dim X_1\;\dots\;\dim X_n}
    \prod_{i=1}^n s_i^{2\dim X_i}(\deg\widebar X_i)\;.\tag\03.6.11.3$$
Here the symbol in parentheses on the right denotes a multinomial coefficient.
By \cite{V~3, Lemma~6.1}, $L_{\epsilon,\bold s}$ is ample.
Hence, if $d$ is sufficiently large then all
the higher cohomology groups vanish, giving
$$h^0(W_{\bold s},dL_{\epsilon,\bold s})
  = \frac{\bigl((dL_{\epsilon,\bold s})^{\sum\dim X_i}\bigr)}{(\sum\dim X_i)!}
    + o\bigl(d^{\sum\dim X_i}\bigr)\;.$$
Combining this with (\03.6.11.3) gives (\03.6.11.1).\qed
\enddemo

\lemma{\03.6.12}
Let $\Gamma(\Cal V,dM_{\epsilon,\bold s})$,
$\Gamma(W_{\bold s},dL_{\epsilon,\bold s})$,
and $\Gamma_\delta\Bigl(\prod\widebar X_i,dM_{\epsilon,\bold s}\Bigr)$
be identified with submodules of
$\Gamma(W_{\bold s},d\pi_{\bold s}^{*}M_{\epsilon,\bold s})$
via $\pi_{\bold s}^{*}$ and (\03.6.8).  Then there is a constant $c>0$,
depending only on $A$, $D$, $L$, $\ell$, $\epsilon$, $\dim X_1,\dots,\dim X_n$,
and the bounds on $\deg\widebar X_1,\dots,\deg\widebar X_n$, such that
if $d$ is sufficiently large and divisible then the rank of the $R$\snug-module
$$\Gamma(\Cal V,dM_{\epsilon,\bold s})
  \cap \Gamma(W_{\bold s},dL_{\epsilon,\bold s})
  \cap \Gamma_\delta\Bigl(\prod\widebar X_i,dM_{\epsilon,\bold s}\Bigr)
\tag\03.6.12.1$$
is bounded from below by $ch^0(W_{\bold s},dL_{\epsilon,\bold s})$.
\endit

\demo{Proof}
Let $Y_i=D\cap\widebar X_i$.  We first show that the upper bound
$$h^0\Bigl(\prod\widebar X_i,dM_{\epsilon,\bold s}\Bigr)
  - h^0_\delta\Bigl(\prod\widebar X_i,dM_{\epsilon,\bold s}\Bigr)
  \le \frac{\delta^n}{n!}\prod_i ds_i^2
    \cdot \prod_i(5ds_i^2)^{\dim Y_i}\cdot\frac{\deg Y_i}{(\dim Y_i)!}
  \tag\03.6.12.2$$
holds.  Let $L'=(4+\epsilon)\rho^{*}L_0+(2+\epsilon)L_1\le5L$.  As noted
in \cite{F}, it suffices to prove the inequality
$$h^0\Bigl(\prod Y_i, \sum ds_i^2\pr_i^{*}L' - \sum e_i\pr_i^{*}D\Bigr)
  \le h^0\Bigl(\prod Y_i, \sum ds_i^2\pr_i^{*}L'\Bigr)$$
for all tuples $(e_1,\dots,e_n)\in\Bbb N^n$
satisfying $e_1/ds_1^2+\dots+e_n/ds_n^2\le\delta$.
(Here $\Bbb N=\{0,1,2,\allowmathbreak\dots\}$.)
This follows from the facts that a translate of $D$ is algebraically equivalent
to an effective divisor on $Y_i$, and that $L'$ is ample.

Then, (\03.6.5), (\03.6.11.1), (\03.6.12.2), and the inequality
$\deg Y_i\le\ell\deg\widebar X_i$ imply that the rank of the module (\03.6.12.1)
is bounded from below by $ch^0(W_{\bold s},dL_{\epsilon,\bold s})$
for some $c>0$.\qed
\enddemo

Next we put three metrics on $\Cal O(M_{\epsilon,\bold s})$.
First, fix a metric on $\Cal O(L_0)$ whose curvature is translation invariant.
For $m=1,\dots,\mu$
let $[0]_m$ and $[\infty]_m$ be the divisors on $\widebar A$ corresponding
to the divisors $\pr_m^{*}[0]$ and $\pr_m^{*}[\infty]$, respectively,
on $(\Bbb P^1)^\mu$.  Then $L_1=\sum_{i=1}^\mu([0]_m+[\infty]_m)$.
By \cite{V~3, Prop.~2.6} at all places $v$ there are
Weil functions $\lambda_{m,v}$ for $[0]_m-[\infty]_m$ satisfying
$$\lambda_{m,v}(P+Q)=\lambda_{m,v}(P)+\lambda_{m,v}(Q)
  \qquad\text{for all $P,Q\in A(\Bbb C_v)$}.$$
(Here $\Bbb C_v$ denotes the completion of the algebraic closure of the
completion $k_v$ of $k$ at $v$; it is algebraically closed.)
These can be used as in \cite{V~3, 2.8} to define smooth metrics
on $\Cal O(L_1)$.
This defines a smooth metric $\|\cdot\|_v$ on $\Cal O(M_{\epsilon,\bold s})$
via the expression (\03.6.7).  Next, for each $v\in S$ let $\|\cdot\|_v'$
be the singular metric constructed in \cite{V~3, 10.5--10.6},
so that $\|\cdot\|'_v$ is equivalent to
the metric on $\Cal O(dL_{\epsilon,\bold s})$ via the embedding (\03.6.8).
Finally, let
$$\|\gamma\|''_v
  = \frac{\|\gamma\|_v}{\sum\exp(-\delta ds_i^2\pr_i^{*}\lambda_{D,v})}\;.
  \tag\03.6.13$$
If a section has index $\ge\delta$ at $D\times\dots\times D$, then this
singular metric remains bounded; one may take this as a definition
of the index at $D\times\dots\times D$.
Note that the singular metrics $\|\cdot\|'$ and $\|\cdot\|''$ are of the
form \cite{V~3, 10.5}; therefore \cite{V~3, Lemma~11.2} applies.

Note also that we take the intersection in (\03.6.12.1) instead of combining
$\|\cdot\|'$ and $\|\cdot\|''$; this is because of problems at infinity
as illustrated in the examples in Sect. \01.  This is why (\03.6.2) is
necessary.

As in \cite{V~3, Thm.~12.4 and Remark~12.6}, we obtain a small section:

\thm{\03.6.14}
For all tuples $\bold s=(s_1,\dots,s_n)$ of positive rational numbers
and for all sufficiently large (and divisible) $d$ (depending on $\bold s$),
there exists a section $\gamma\in\Gamma(\Cal V,dM_{\epsilon,\bold s})$
such that $\|\gamma\|'$ and $\|\gamma\|''$ are bounded
and such that the inequality
$$\prod_{v\mid\infty}\|\gamma\|_{\sup,v}
  \le \exp\Biggl(cd\sum_{i=1}^n s_i^2\Biggr)\tag\03.6.14.1$$
holds.  Here the constant $c$ is independent of $\bold s$ and $d$.
\endit

\demo{Proof}
This proof is a matter of obtaining bounds for volumes
of various lattices in the diagram
$$\CD 0 & \to & \Gamma(\Cal V, dM_{\epsilon,\bold s})
    @>>> \Gamma\Bigl(\Cal V, d\sum s_i^2\pr_i^{*}L'\Bigr)^a
    @>>> \Gamma\Bigl(\Cal V, d\sum s_i^2\pr_i^{*}L''\Bigr)^b \\
  && \cup && \cup \\
  && \Gamma'(\Cal V, dM_{\epsilon,\bold s})
    @>>> \Gamma'\Bigl(\Cal V, d\sum s_i^2\pr_i^{*}L'\Bigr)^a \\
  && \cup && \cup \\
  && \Gamma''(\Cal V, dM_{\epsilon,\bold s})
    @>>> \Gamma''\Bigl(\Cal V, d\sum s_i^2\pr_i^{*}L'\Bigr)^a\rlap{\quad.}
  \endCD$$
This is as in \cite{V~3}:  the top row is the Faltings complex, with
$$L' = 4\rho^{*}L_0 + 2L_1 + \epsilon L$$
and
$$L'' = 8\rho^{*}L_0 + 2L_1 + \epsilon L\;;$$
the symbols $\Gamma'$ in the middle row denote the submodules
of sections $\gamma$ for which $\|\gamma\|'$ is bounded; and
the symbols $\Gamma''$ in the bottom row denote the submodules
of sections $\gamma$ for which both $\|\gamma\|'$ and $\|\gamma\|''$
are bounded.

The proof is exactly the same as in \cite{V~3}, except that the
fifth paragraph is repeated because of the extra row in the above diagram.\qed
\enddemo

By (\03.6.13), \cite{V~3, Prop.~10.10}, (\03.6.10.4), (\03.6.1),
and the choice of the $s_i$, we have
$$\split \|\gamma(P_1,\dots,P_n)\|_w
  &= \|\gamma(P_1,\dots,P_n)\|''_w
    \sum_{i=1}^n\exp(-\delta ds_i^2\lambda_{D,w}(P_i)) \\
  &\le \|\gamma\|''_{\sup,w}
    \cdot n\max_{1\le i\le n}\exp(-\delta ds_i^2\lambda_{D,w}(P_i)) \\
  &\le \|\gamma\|_{\sup,w}\cdot\exp(-\delta d\kappa)
    \cdot\exp\Biggl(cd\sum_{i=1}^n s_i^2\Biggr)
  \endsplit$$
and therefore
$$-\log\|\gamma(P_1,\dots,P_n)\|_w
  \ge -\log\|\gamma\|_{\sup,w} + d\kappa\delta - cd\sum_{i=1}^n s_i^2\;.
  \tag\03.6.15$$
Likewise, letting $\alpha_{vmi}=\exp(-\pr_i^{*}\lambda_{m,v})$ and applying
\cite{V~3, (10.6) and Prop.~10.10} to $\|\cdot\|'$ gives
$$\split &-\log\|\gamma(P_1,\dots,P_n)\|_v \\
  &\qquad\qquad\ge -\log\|\gamma\|_{\sup,v}
      + \sum_{m=1}^\mu\sum_{i=2}^n -\log
        \frac{\Bigl(\alpha_{vm,i-1}^{ds_{i-1}^2}+\alpha_{vmi}^{ds_i^2}\Bigr)^2}
        {\Bigl(1+\alpha_{vm,i-1}^{ds_{i-1}^2}\Bigr)^2
          \Bigl(1+\alpha_{vmi}^{ds_i^2}\Bigr)^2}
    - cd\sum_{i=1}^n s_i^2\endsplit\tag\03.6.16$$
for all $v\in S$; cf\. \cite{V~3, Prop.~12.5}.  If the $P_i$ are
chosen so that \cite{V~3, 13.5} holds, then by \cite{V~3, 13.7} we have
$$\split &ds_1^2h_{L_1}(P_1) + 2\sum_{i=2}^{n-1}ds_i^2h_{L_1}(P_i)
    + ds_n^2h_{L_1}(P_n) \\
  &\qquad\le\frac1{[k:\Bbb Q]}\sum_{v\in S}\sum_{m=1}^\mu\sum_{i=2}^n -\log
        \frac{\Bigl(\alpha_{vm,i-1}^{ds_{i-1}^2}+\alpha_{vmi}^{ds_i^2}\Bigr)^2}
        {\Bigl(1+\alpha_{vm,i-1}^{ds_{i-1}^2}\Bigr)^2
          \Bigl(1+\alpha_{vmi}^{ds_i^2}\Bigr)^2} \\
  &\qquad\qquad+ d(n-1)\epsilon_1 + cd\sum_{i=1}^n s_i^2\;.
  \endsplit\tag\03.6.17$$
But also, choosing $P_1,\dots,P_n$ so that (\03.6.2) holds for sufficiently
small $\eta$, we have
$$\frac1{[k:\Bbb Q]}\sum_{m=1}^\mu\sum_{i=2}^n -\log
        \frac{\Bigl(\alpha_{wm,i-1}^{ds_{i-1}^2}+\alpha_{wmi}^{ds_i^2}\Bigr)^2}
        {\Bigl(1+\alpha_{wm,i-1}^{ds_{i-1}^2}\Bigr)^2
          \Bigl(1+\alpha_{wmi}^{ds_i^2}\Bigr)^2}
  \le d\epsilon_2\tag\03.6.18$$
for some $\epsilon_2>0$ depending only on $\eta$, $\mu$, $n$, and $[k:\Bbb Q]$.
Adding (\03.6.16) for $v\in S\setminus\{w\}$ to (\03.6.15)
and applying (\03.6.14.1), (\03.6.18), and (\03.6.17) then gives
$$\split &\frac1{[k:\Bbb Q]}\sum_{v\in S}-\log\|\gamma(P_1,\dots,P_n)\|_v \\
  &\qquad\ge \frac1{[k:\Bbb Q]}\left(
      \sum_{v\in S\setminus\{w\}}\sum_{m=1}^\mu\sum_{i=2}^n -\log
        \frac{\Bigl(\alpha_{wm,i-1}^{ds_{i-1}^2}+\alpha_{wmi}^{ds_i^2}\Bigr)^2}
        {\Bigl(1+\alpha_{wm,i-1}^{ds_{i-1}^2}\Bigr)^2
          \Bigl(1+\alpha_{wmi}^{ds_i^2}\Bigr)^2}
    + d\kappa\delta\right) \\
  &\qquad\qquad+ \frac1{[k:\Bbb Q]}\sum_{v\in S}-\log\|\gamma\|_{\sup,v}
    - cd\sum_{i=1}^n s_i^2 \\
  &\qquad\ge \frac1{[k:\Bbb Q]}\left(\sum_{v\in S}\sum_{m=1}^\mu\sum_{i=2}^n
        -\log
        \frac{\Bigl(\alpha_{wm,i-1}^{ds_{i-1}^2}+\alpha_{wmi}^{ds_i^2}\Bigr)^2}
        {\Bigl(1+\alpha_{wm,i-1}^{ds_{i-1}^2}\Bigr)^2
          \Bigl(1+\alpha_{wmi}^{ds_i^2}\Bigr)^2}
    + d\kappa\delta\right) \\
  &\qquad\qquad- d\epsilon_2 - cd\sum_{i=1}^n s_i^2 \\
  &\qquad\ge d\Biggl(\frac{\kappa\delta}{[k:\Bbb Q]} - (n-1)\epsilon_1
      - \epsilon_2 + s_1^2h_{L_1}(P_1)
      + 2\sum_{i=2}^{n-1}s_i^2h_{L_1}(P_i) + s_n^2h_{L_1}(P_n)\Biggr) \\
  &\qquad\qquad - cd\sum s_i^2\;.\endsplit$$
On the other hand, as in \cite{V~3, 13.6}, for suitably chosen $P_1,\dots,P_n$,
we have
$$\frac1{[k:\Bbb Q]}\deg M_{\epsilon,\bold s}\restrictedto E
  \le (n-1)\epsilon_1+n\epsilon
    + s_1^2h_{L_1}(P_1) + 2\sum_{i=2}^{n-1}s_i^2h_{L_1}(P_i) + s_n^2h_{L_1}(P_n)
    + c\sum s_i^2$$
on the arithmetic curve $E$ corresponding to $(P_1,\dots,P_n)$.

Combining these two inequalities gives
$$\frac1{[k:\Bbb Q]}\sum_{v\notin S}-\log\bigl\|\gamma\restrictedto E\bigr\|_v
  \le d\Bigl(n\epsilon-\frac{\kappa\delta}{[k:\Bbb Q]}
    + 2(n-1)\epsilon_1 + \epsilon_2\Bigr) + cd\sum s_i^2\;.$$
By (\03.6.4) this gives a negative upper bound if $\epsilon_1$ and $\epsilon_2$
are sufficiently small, leading to a positive lower bound for
the index of $\gamma$ at $(P_1,\dots,P_n)$
(with multiplicities $ds_1^2,\dots,ds_n^2$).
The argument then concludes by applying the product theorem in the usual
way, as in \cite{V~2, Sect.~18}.\qed
\enddemo

\beginsection{\04}{Proof of Theorem \00.2}

This section uses the notation given in the introduction of the paper.

First, we may assume that the Ueno fibration is trivial.
This is because the theorem is preserved under pulling back by quotient
morphisms.  We may enlarge $k$ so that the toric part of $A$ splits;
\ie, the exact sequence (\00.1) holds already over $k$.  Finally, it will
suffice to assume that $\Cal V(R_S)$ is Zariski-dense and obtain
a contradiction.  To see this, apply Theorem \00.3 to any irreducible
component of the closure of $\Cal V(R_S)$ and then proceed by induction
on the dimension.

Let $\widebar A$, $\pi\:\widetilde A\to\widebar A$, $\widebar D$,
and $\widetilde D$
be as in Theorem \02.4.  For $v\in S$ let $\lambda_{\widetilde D,v}$
be a Weil function for $\widetilde D$ on $\widetilde A$.
Then there exist $w\in S$ and $\kappa_0>0$ such that
$$\lambda_{\widetilde D,w}(\widetilde P)
  \ge \kappa_0 h_{\widetilde D}(\widetilde P)$$
for all $\widetilde P\in\widetilde A$ corresponding to elements $P\in\Cal S_0$,
where $\Cal S_0\subseteq\Cal V(R_S)$ is Zariski-dense in $A$.  We now claim
that, after shrinking to a possibly smaller $\Cal S\subseteq\Cal S_0$,
the inequality
$$\lambda_{\widetilde D,w}(\widetilde P) \ge \kappa h_{\widebar D}(P)\tag\04.1$$
holds as well.  If not, then by (\02.4.1)
$h_F(\widetilde P)\ge\kappa h_{\widebar D}(P)$ for some $\kappa>0$.
But since $\pi(\Supp F)$ is contained in $\widebar{\Supp D}$,
we have $h_F(\widetilde P)=\sum_{v\in S}\lambda_{F,v}(\widetilde P)+O(1)$
and therefore
$$\lambda_{F,w}(\widetilde P)\ge\kappa h_{\widebar D}(P)$$
for some $\kappa>0$ and some $w\in S$ (after shrinking $\Cal S$).
Pushing down to $A_0$ then gives an inequality which contradicts
Theorem \03.6 in the case where $\mu=0$ (which is also \cite{F, Thm.~2}).
Thus (\04.1) holds.

For $A$\snug-orbits $T\ne A$ of $\widebar A$ let $\widebar T$ be the closure and
let $\lambda_{\widebar T,w}$ be some logarithmic distance function for
$\widebar T$.  Choose an orbit $T$ of minimal dimension such that
$$\lambda_{\widebar T,w}(P) \ge \eta h_{\widebar D}(P)$$
for some $\eta>0$ and for all $P$ in some Zariski-dense subset $\Cal S'$
of $\Cal S$.  If there is no such orbit then let $T=A$.

Let $p\:A\to T$ be the restriction to $A$ of the equivariant projection
defined in Proposition \02.5.  Let $\widetilde T=\pi^{-1}(\widebar T)$;
then $p$ lifts to a morphism $\tilde p\:\pi^{-1}(A)\to\widetilde T$.
Points $\tilde p(\widetilde P)$ for $\widetilde P\in\widetilde A$
lying over $P\in\Cal S'$ approach $\widetilde D\cap\widetilde T$
as in (\03.6.1):  for a suitable $\kappa'>0$, we have
$$\lambda_{\widetilde D,w}(\tilde p(\widetilde P))
  \ge \kappa' h_{\widebar D}(P)\;.\tag\04.2$$
But also, by Proposition \02.3, there is a constant $c_1>0$ such that
$$h_{\widebar D}(P)\ge c_1 h_{\widebar D}(p(P))+O(1)$$
for all $P\in\Cal S$.  We may therefore replace $h_{\widebar D}(P)$ in (\04.2)
with $h_{\widebar D}(p(P))$.  This gives (\03.6.1) for $p(\Cal S')$ on $T$,
since $\lambda_{\widetilde D,w}$ is a generalized Weil function on $T$.
The condition (\03.6.2) also holds, by minimality of the choice of $T$,
and by (\03.5).  This leads to a contradiction, by Theorem \03.6 applied
to $\widebar T$.\qed

\beginsection{\05}{Some additional geometry}

The similarity between the conclusions of Theorems \00.2 and \00.3 suggests
that some of the results traditionally proved for closed subvarieties of
semiabelian varieties could be proved for closed subvarieties minus divisors,
too.  This section generalizes results of Ueno and Fujita on the logarithmic
Kodaira dimension of such varieties.

\smallskip\bgroup\narrower\it Unless otherwise specified, all varieties are
over an algebraically closed field $k$ of characteristic zero.
\par\egroup\smallbreak

Many of these results probably extend to positive characteristic, but
additional work would be needed due to the unavailability of resolution of
singularities in positive characteristic.

For general references on group varieties, see Rosenlicht \cite{R};
for other references on closed subvarieties of semiabelian varieties, see
Abramovich \cite{A}.

\thm{\05.1}  Let $k$ be any field, let $X$ be a complete nonsingular variety
over $k$, and let $D$ be a divisor on $X$.
Suppose that $H^0(X,mD)\ne0$ for some $m\in\Bbb Z_{>0}$.  Let $m_0$
be the index of the subgroup of $\Bbb Z$ generated by all such $m$.
Then there exist constants $c_2\ge c_1>0$ and an integer $\kappa$ with
$0\le\kappa\le\dim X$ such that
$$c_1 m^\kappa \le h^0(X,mm_0D) \le c_2 m^\kappa$$
for all sufficiently large $m$.
\endit

\demo{Proof}  See \cite{I~3, Thm.~10.2}.\qed
\enddemo

\defn{\05.2}  Let $k$ be any field and let $X$ be a complete nonsingular
variety over $k$.
\roster
\myitem a.  Let $D$ be a divisor on $X$.  The {\bc $D$\snug-dimension} of $X$,
denoted $\kappa(X,D)$, is the number $\kappa$ in Theorem \05.1 if $D$
satisfies the conditions of that theorem; otherwise it is $-\infty$.
\myitem b.  Let $\Cal L$ be a line sheaf on $X$.  Then the
{\bc $\Cal L$\snug-dimension} of $X$, denoted $\kappa(X,\Cal L)$, is defined
to be $\kappa(X,D)$ for any divisor $D$ on $X$ such that $\Cal L\cong\Cal O(D)$.
\myitem c.  A divisor $D$ (resp\. line sheaf $\Cal L$) on $X$ is {\bc big}
if $\kappa(X,D)$ (resp\. $\kappa(X,\Cal L)$) equals $\dim X$.
\endroster
\endit

A divisor $D$ on $X$ is big if and only if $h^0(X,mD)\gg m^{\dim X}$
for all sufficiently large and divisible integers $m$.

\defn{\05.3}  Let $V$ be a nonsingular quasi-projective variety.
\roster
\myitem a.  A {\bc smooth completion of $V$} is an open
immersion $V\hookrightarrow X$ into a nonsingular projective variety $X$
such that $D:=X\setminus V$ is a normal crossings divisor (taken with all
multiplicities equal to one).
Such a smooth completion will often be denoted by a pair $(X,D)$.
\myitem b.  The {\bc logarithmic canonical divisor} of $V$ on $X$ is
$K_{(X,D)}:=K_X+D$.  If $(X,D)$ is clear from the context, then it may
be denoted $K_V$ and called simply the logarithmic canonical divisor
of $V$.
\myitem c.  The {\bc logarithmic Kodaira dimension} of $V$ is the number
$\kappa(X,K_{(X,D)})$; by \cite{I~3, Thm.~11.2} it depends only on $V$.
It is denoted $\bar\kappa(V)$.
\myitem d.  We say that $V$ is of {\bc logarithmic general type}
if $\bar\kappa(V)=\dim V$; \ie, if $K_V$ is big.
\myitem e.  Let $\pi\:V\to W$ be a morphism to a nonsingular quasi-projective
variety $W$, let $V\hookrightarrow X$ and $W\hookrightarrow Y$ be smooth
completions such that $\pi$ extends to a morphism $\bar\pi\:X\to Y$,
and let $K_V$ and $K_W$ be the logarithmic canonical divisors of $V$ and $W$,
respectively.  Then the {\bc relative logarithmic canonical divisor}
of $V$ over $W$ with respect to $\bar\pi$ is
$K_{V/W}:=K_V-\bar\pi^{*}K_W$.  Again, mention of $\bar\pi$ may be
omitted if it is clear from the context.
\endroster
\endit

The first result of this section is that if a closed subvariety of a
semiabelian variety, minus a divisor, has trivial Ueno fibration,
then it is of logarithmic general type.  For closed
subvarieties this was proved by Iitaka \cite{I~1} and \cite{I~2}; the proof
here is an easy adaptation of that proof.  Since Iitaka's exposition often
leaves out details, however, we provide a more complete proof here.  This
also provides the opportunity to change the proof a little, by replacing
a cardinality argument with an argument on the field of definition of a
group subvariety.

The first step in the proof consists of proving it in the special case
when the closed subvariety is the whole semiabelian variety; in other words,
a semiabelian minus a divisor $D$ with $B(D)=0$ has logarithmic general
type.  In the end of the paper \cite{I~2}, Iitaka remarks that this
was proved by T. Fujita, but I have been unable to find a reference.
Therefore, we give a proof here, adapting a proof of Mumford
\cite{Mu, \S6, pp\. 60--61} for part of the way.

\lemma{\05.4}
\roster
\myitem a.  Let $\widebar G$ be a toric smooth completion of $\Gm^\mu$.
Then the logarithmic canonical divisor of $\Gm^\mu$ on $\widebar G$
is trivial, and the isomorphism between the canonical line sheaf and the
trivial line sheaf on $\widebar G$ is translation invariant.
\myitem b.  Let $D$ be an effective divisor on $\Gm^\mu$ with trivial
Ueno fibration.  Then there exists a toric smooth completion $\widebar G$
of $\Gm^\mu$, such that the closure $\widebar D$ of $D$ in $\widebar G$
does not contain any $\Gm^\mu$\snug-orbit of $\widebar G$.
\endroster
\endit

\demo{Proof}
A nonsingular toric variety $\widebar G$ with principal orbit $\Gm^\mu$
can be described by giving an open cover in which each open subset is
isomorphic to $\Bbb A^\mu$, and the action of $\Gm^\mu$ is by monomials.
Therefore, the divisor $\widebar G\setminus\Gm^\mu$ is a normal crossings
divisor (and hence $\widebar G$ is a smooth completion of $\Gm^\mu$).
Also, the differential $\mu$\snug-form $dx_1/x_1\wedge\dots\wedge dx_\mu/x_\mu$
is a generator on this open subset for the canonical line sheaf of $\Gm^\mu$
on $\widebar G$.  This generator is translation invariant, and is a nonzero
constant multiple of the corresponding generator over any other such open
subset.  This proves part (a).

For the proof of (b) we first recall some definitions from the theory of
toric varieties.  A {\bc fan} is a finite set $\Sigma$ of polyhedral cones
in $\Bbb R^\mu$ such that every $\sigma\in\Sigma$ is a closed, rational,
polyhedral cone not containing any nontrivial linear subspace; all faces
of all $\sigma\in\Sigma$ also lie in $\Sigma$; and for all
$\sigma,\sigma'\in\Sigma$, $\sigma\cap\sigma'$ is a face of $\sigma$ and
of $\sigma'$.
A {\bc barycentric subdivision} of $\Sigma$ associated to a rational ray
$\lambda$ is the fan $\Sigma'$ consisting of all cones in $\Sigma$ not
containing $\lambda$, plus the convex hulls of $\lambda$ with each face
of each $\sigma\in\Sigma$ containing $\lambda$.

By Lemma \02.4.2 there exists an equivariant completion $G_1$ of $\Gm^\mu$
such that the closure $D_1$ of $D$ in $G_1$ is Cartier and ample.
This is a toric variety.  Moreover, $D_1$ does not contain any orbit
of $G_1$.

By \cite{D, 8.1}, $G_1$ can be desingularized by applying a finite
sequence of barycentric subdivisions to the corresponding fan $\Sigma$.
Consider one such barycentric subdivision:  let $\lambda$ be a rational ray,
and let $\phi\:X'\to X$ be the corresponding morphism of toric varieties.
Orbits in $X'$ map onto orbits in $X$, so if the closure $\widebar D_X$
of $D$ in $X$ does not contain any orbit in its support, then the same is
true of the closure $\widebar D_{X'}$ of $D$ in $X'$.
In particular, this applies to the exceptional set of $\phi$,
which is the closure of the orbit corresponding to the ray $\lambda$.
Thus $\widebar D_X$ and $\widebar D_{X'}$ are related by pull-back of
Cartier divisors:  $\widebar D_{X'}=\phi^{*}\widebar D_X$. 

Let $\phi\:\widebar G\to G_1$ be the desingularization corresponding to
the composite of these barycentric subdivisions.  By induction, the closure
$\widebar D$ of $D$ in $\widebar G$ does not contain any orbit under the
action of $\Gm^\mu$.  Moreover, $\widebar G$ is a toric nonsingular
equivariant completion of $\Gm^\mu$.  Thus part (b) holds.\qed
\enddemo

\cor{\05.5}  Let $\rho\:X\to Y$ be a fiber bundle with fiber $\Gm^\mu$,
where $X$ and $Y$ are varieties,
and let $D$ be an effective divisor on $X$ with $B(\Gm^\mu,D)=0$.
Then there exists a toric smooth completion $\widebar G$ of $\Gm^\mu$
such that the closure $\widebar D$ of $D$ in the corresponding (relative)
completion $\widebar X$ of $X$ does not contain any subset corresponding
to a $\Gm^\mu$\snug-orbit of $\widebar G$.
\endit

\demo{Proof}  Applying Lemma \05.4b to the generic fiber of $\rho$ gives
a toric smooth completion $\widebar G_{\eta}$ of $\Gm^\mu$ over the
generic point $\eta$ of $Y$.  Since toric varieties are defined by
discrete data, $\widebar G_{\eta}$ is of the form $\widebar G\times_k\eta$
for some toric smooth completion $\widebar G$ of $\Gm^\mu$ over $k$.
This is the desired $\widebar G$.\qed
\enddemo

\lemma{\05.6}  Let $\pi\:X\to Y$ be a fiber bundle with fiber $B$,
where $X$ and $Y$ are nonsingular varieties and $B$ is a semiabelian variety.
Assume that there exists a smooth completion $(\widebar Y,E)$ of $Y$
such that $\pi$ extends to a fiber bundle over $\widebar Y$.
Then there exists a smooth completion $(\widebar X,D)$ of $X$ such that
$\pi$ extends to a fiber bundle $\bar\pi\:\widebar X\to\widebar Y$,
and such that $K_{\widebar X}+D=\bar\pi^{*}(K_{\widebar Y}+E)$.
Consequently, $\bar\kappa(X)=\bar\kappa(Y)$.  Moreover, $\bar\pi$ can
be constructed so that its fiber is an equivariant completion of $B$
corresponding to a toric smooth completion of the toric part of $B$.
\endit

\demo{Proof}  We first claim that if $A$ is a semiabelian variety,
if $\widebar A$ is an equivariant completion as in Lemma \02.2
with $\widebar G$ a toric smooth completion of $\Gm^\mu$, and if
$D=\widebar A\setminus A$, then $K_{(\widebar A,D)}$ is trivial, and
the isomorphism $\Cal O(K_{(\widebar A,D)})\cong\Cal O_{\widebar A}$
commutes with translation.
If $A$ is an abelian variety, then this is classical, since in fact
$\Omega^1_{A/k}\cong\Cal O_A^{\dim A}$ and that isomorphism commutes with
translation.  If $A=\Gm^\mu$, then this follows from Lemma \05.4a.
In general, we have $\rho^{-1}(U)\cong U\times\Gm^\mu$,
so $K_{\rho^{-1}(U)/U}=0$ for open $U\subseteq A_0$ in an open covering of
$A_0$, and the transition functions between these isomorphisms consist
of translations on $\Gm^\mu$, so the isomorphisms $K_{\rho^{-1}(U)/U}=0$
patch together to give us $K_{A/A_0}=0$.  This isomorphism is invariant
under translations by $A$, since the same is true on suitable open subsets
of the sets $\rho^{-1}(U)$.

To prove the lemma itself, let $(\widebar Y,E)$ be as assumed, and let
$\widebar B$ be an equivariant completion of $B$ as above.
This determines a smooth $B$\snug-equivariant completion $(\widebar X,D)$
of $X$, as in Lemma \02.2.  As before, a patching argument then gives
$K_{\widebar X}+D=\bar\pi^{*}(K_{\widebar Y}+E)$.
The assertion $\bar\kappa(X)=\bar\kappa(Y)$ follows immediately from this.\qed
\enddemo

\lemma{\05.7}  Let $A$ be a semiabelian variety and let $B$ be a semiabelian
subvariety.  Then there exist equivariant completions of $A$ and $A/B$
such that the canonical map $A\to A/B$ extends to a morphism between the
completions.
\endit

\demo{Proof}  As usual, let $\rho\:A\to A_0$ be the maximal abelian
quotient of $A_0$.  Choose an isomorphism $\Ker\rho\cong\Gm^\mu$ such that
the first $\mu(B)$ factors correspond to the maximal torus in $B$.
One can then use the completion $\widebar A$ corresponding to the completion
$\Gm^\mu\hookrightarrow(\Bbb P^1)^\mu$.\qed
\enddemo

\lemma{\05.8} (Theorem of the Square)  Let $A$ be a semiabelian variety
and let $\widebar A$ be an equivariant completion of $A$ corresponding to
a nonsingular toric equivariant completion of $\Gm^\mu$ (in the notation of
(\00.4)).  For $x\in A$ let $T_x\:\widebar A\to\widebar A$ denote
translation by $x$.  Then for all line sheaves $\Cal L$ on $\widebar A$
and all $x,y\in A(k)$,
$$\Cal L \otimes T_{x+y}^{*}\Cal L
  \cong T_x^{*}\Cal L \otimes T_y^{*}\Cal L\;.$$
\endit

\demo{Proof}  Let $\widebar G$ be the equivariant completion
of $\Gm^\mu$ mentioned above.  Then (since all divisors on $\Gm^\mu$
are principal) every divisor on $\widebar G$ is linearly equivalent to
a divisor supported on closures of orbits.
Therefore $\Cal L\cong\Cal O(D_1+\rho^{*}D_2)$, where $D_1$ is a divisor
supported only on subsets of $\widebar A$ corresponding to orbits of
$\widebar G$, and $D_2$ is a divisor on $A_0$.  The lemma then follows,
since $D_1$ is invariant under translation, and since the theorem of the
square holds on $A_0$.\qed
\enddemo

The following is an adaptation of a result proved by Mumford for abelian
varieties; cf\. \cite{Mu, \S6, pp\. 60--61}.

\lemma{\05.9}  Let $A$ be a semiabelian variety, let $D$ be an effective
divisor on $A$, and let $\widebar A$ be an equivariant completion of $A$
corresponding to a completion $\widebar G$ of $\Gm^\mu$ satisfying the
conditions of Corollary \05.5.  Let $\widebar D$ be the closure of $D$
in $\widebar A$.  Then
\roster
\myitem a.  the linear system $|2\widebar D|$ is base-point free;
\myitem b.  if $B(D)=0$, then the morphism $\widebar A\to\Bbb P^N$ induced
by $|2\widebar D|$ is generically finite; and
\myitem c.  if $B(D)=0$, then $\widebar D$ is big.
\endroster
\endit

\demo{Proof}  Lemma \05.8 implies that
$T_x^{*}\widebar D+T_{-x}^{*}\widebar D\sim 2\widebar D$ for all $x\in A(k)$,
so given any $P\in\widebar A$ it suffices to find $x\in A(k)$ such that
$P\pm x\notin\Supp\widebar D$.  For suitably generic choices of $x$,
$\Supp\widebar D$ does not contain any orbit of $\Gm^\mu$ in the fiber of
$\rho$ containing $P+x$ or $P-x$.  This condition depends only on $\rho(x)$.
The condition $P\pm x\notin\Supp\widebar D$ is then satisfied for a
generic choice of $x$ within such a fiber.  This proves (a).

Now suppose $B(D)=0$, and let $\phi$ be a morphism $\widebar A\to\Bbb P^N$
induced by $|2\widebar D|$.  If $\phi$ is not generically finite, then
there is an integral curve $C$ in $\widebar A$, meeting $A$, such that
$\phi(C)$ is a point.
Then for all $x$, $\Supp\bigl(T_x^{*}\widebar D+T_{-x}^{*}\widebar D\bigr)$
either contains $C$, or is disjoint from $C$.  In particular, it is
disjoint from $C$ for almost all $x$.  Hence the same is true of
$\Supp T_x^{*}\widebar D$, and of $T_x^{*}D_0$ for all irreducible
components $D_0$ of $\widebar D$.

We now claim that all such components $D_0$ are invariant under translation
by $x_2-x_1$ for all $x_1,x_2\in C\cap A$.  Indeed, since all divisors
$T_x^{*}D_0$ are algebraically equivalent, their restrictions to $C$ must
have the same degree, which must be zero since $C$ is usually disjoint from
$T_x^{*}D_0$.  Let $x_1,x_2\in C\cap A$ and $y\in D_0$.
Then $x_1\in T_{y-x_1}^{*}D_0$, so also $x_2\in T_{y-x_1}^{*}D_0$,
and therefore $y\in T_{x_2-x_1}^{*}D_0$.  This holds for all $y\in D_0$,
so $D_0\subseteq T_{x_2-x_1}^{*}D_0$.  By symmetry they are equal, thus
proving the claim.

But now it follows that $\widebar D$ is invariant under translation by
$x_2-x_1$ for all $x_1,x_2\in C\cap A$, contradicting the assumption that
$B(D)=0$.  Thus $\phi$ is generically finite.

Finally, $\widebar D$ is big, because $2\widebar D$ is the pull-back of
a hyperplane via the generically finite morphism $\phi$.\qed
\enddemo

\lemma{\05.10}  Let $V$ be a nonsingular quasi-projective variety and
let $X$ be a smooth completion of $V$.  Let $D$ be an effective
divisor on $X$, all of whose irreducible components meet $V$.
Assume also that $\Supp D$ does not contain any irreducible local
intersection of components of $X\setminus V$.
Then there exists a smooth completion $X'$ of $V\setminus D$ admitting
a morphism $\pi\:X'\to X$ such that the relative logarithmic canonical
divisor of $\pi^{-1}(V\setminus D)$ over $V$ is effective, with support
equal to $\pi^{-1}(\Supp D)$.
\endit

\demo{Proof}
Let $F$ be the divisor $X\setminus V$ (with all multiplicities equal to one).
By Hironaka's resolution of singularities
\cite{Hi, pp\. 142--143, Main Theorem II}
there exists a sequence $X_r\to X_{r-1}\to\dots\to X_0=X$ of blowings-up
such that
\roster
\myitem 1.  $\pi_{i+1}\:X_{i+1}\to X_i$ is the blowing-up of an irreducible
nonsingular subvariety $C_i$ which has normal crossings with $E_i\cup F_i$,
where $F_i$ is the inverse image of $F$ in $X_i$, and $E_i$ is the exceptional
set of the morphism $\pi_1\circ\dots\circ\pi_i$;
\myitem 2.  $C_i$ is contained in the strict transform of $D$ for all $i$;
and
\myitem 3.  The strict transform of $D$ in $X_r$ is nonsingular and has normal
crossings with $E_r\cup F_r$.
\endroster

We claim that for all $i$ the relative logarithmic canonical divisor
of $X_i\setminus(E_i\cup F_i)$ over $X\setminus F$ is
an effective divisor whose support equals $E_i$.  This will be proved by
induction.  It is trivial if $i=0$.  Assume it is true for $i$.
We may assume that $C_i$ is not a divisor.  If $C_i$ is not locally an
intersection of components of $E_i\cup F_i$, then the relative logarithmic
canonical divisor of $X_{i+1}\setminus(E_{i+1}\cup F_{i+1})$ over
$X_i\setminus(E_i\cup F_i)$ is effective, with support equal to the
exceptional divisor of $\pi_{i+1}$, so the inductive hypothesis is true
for $i+1$.  If, on the other hand, $C_i$ is locally an intersection
of components of $E_i\cup F_i$, then at least one of these local components
must be a component of $E_i$, for otherwise $\pi_1\circ\dots\circ\pi_i$
would be \'etale in a neighborhood of the generic point of $C_i$,
and $(\pi_1\circ\dots\circ\pi_i)(C_i)$ would be a local intersection of
components of $F$, contradicting the assumption on $\Supp D$.  Thus,
in this case, the relative logarithmic canonical divisor of 
$X_{i+1}\setminus(E_{i+1}\cup F_{i+1})$ over $X_i\setminus(E_i\cup F_i)$
is zero, and $E_{i+1}=\pi_{i+1}^{-1}(E_i)$, so the inductive hypothesis
is again satisfied.

In particular, this claim holds for $i=r$.  Let $\pi\:X_r\to X$ be the
composition $\pi_1\circ\dots\circ\pi_r$.  Then the relative logarithmic
canonical divisor of
$$X_r\setminus(E_r\cup F_r\cup\pi^{-1}(\Supp D)) = \pi^{-1}(V\setminus D)$$
over $V$ is an effective divisor whose support equals $\pi^{-1}(\Supp D)$.\qed
\enddemo

Finally, we can now prove Fujita's result.

\lemma{\05.11}  Let $A$ be a semiabelian variety, and let $D$ be an effective
divisor on $A$ with $B(D)=0$.  Then $A\setminus D$ is of logarithmic
general type.
\endit

\demo{Proof}  By Lemmas \05.4 and \05.9, there exists a toric smooth
completion $\widebar A$ of $A$ such that the closure $\widebar D$ of $D$
is big, and such that $\widebar D$ does not contain any local intersection
of components of $\widebar A\setminus A$.
By Lemma \05.10, there exists a proper birational morphism
$\pi\:\widetilde A\to\widebar A$ such that $\widetilde A$ is a smooth
completion of $A\setminus D$ and such that the relative logarithmic canonical
divisor of $\pi^{-1}(A\setminus D)$ (in $\widetilde A$)
over $A$ is effective, with support equal to $\pi^{-1}(\Supp\widebar D)$.
In particular, it is big.  But, by Lemma \05.6 with $Y$ equal to a point,
the logarithmic canonical divisor of $A$ (in $\widebar A$) is zero;
hence the logarithmic canonical divisor of $A\setminus D$ (on $\widetilde A$)
is big.  Thus $A\setminus D$ is of logarithmic general type.\qed
\enddemo

\cor{\05.12} (Fujita)  Let $D$ be an effective divisor on a semiabelian
variety $A$.  Then $\bar\kappa(A\setminus D)\ge0$, with equality if and
only if $D=0$.
\endit

\demo{Proof}  Let $A'=A/B(D)$ and $D'=D/B(D)$.  By Lemma \05.11,
$\bar\kappa(A'\setminus D')=\dim A'$.  By Lemmas \05.6 and \05.7,
$\bar\kappa(A\setminus D)=\bar\kappa(A'\setminus D')$.
Thus $\bar\kappa(A\setminus D)\ge0$, with equality if and only if $B(D)=A$.
The latter condition holds if and only if $D=0$.\qed
\enddemo

\lemma{\05.13}  Let $k$ be any field (of any characteristic, not necessarily
algebraically closed).  Let $A$ be a semiabelian variety over $k$,
with maximal abelian quotient $\rho\:A\to A_0$.
Let $S$ be a geometrically integral scheme over $k$ with $S(k)\ne\emptyset$,
and let $\Cal B$ be a reduced closed subscheme of $A\times_k S$.
Assume that $\Cal B$ is a group subscheme of the group scheme $A\times_k S$
over $S$, that $\Cal B$ is smooth over $S$, that the restriction
of the map $\rho_{(S)}\:A\times_k S\to A_0\times_k S$ to $\Cal B$ is smooth,
and that $\Cal B$ is geometrically connected over $k$.
Then there is a group subvariety $B$ of $A$ such that $\Cal B=B\times_k S$.
(Cf\. \cite{Mi, Prop.~20.3}.)
\endit

\demo{Proof}  Pick $s\in S(k)$, and let $B=\Cal B_s$.  Since $\Cal B$
is reduced, it is sufficient
to show that $\Cal B\times_k\bar k=B\times_k S\times_k\bar k$,
set-theoretically.  Hence we may assume that $k$ is algebraically closed.

Let $K$ be the function field of $S$.  The image of
$\rho_{(K)}\restrictedto{\Cal B_K}\:\Cal B_K\to A_0\times_k K$ is a
connected group subvariety; hence an abelian subvariety; by \cite{Mi, Cor.~20.4}
it is of the form $B_0\times_k K$ for some abelian subvariety $B_0$ of $A_0$.
By smoothness and dimensionality considerations, $\rho_{(S)}$ maps $\Cal B$
onto $B_0\times_k S$.  By shrinking $A$, we may therefore assume that $B_0=A_0$.

Now consider the group $C_K:=\Ker\rho_{(S)}\cap\Cal B_K$.  It is a subgroup
of an algebraic torus; hence by \cite{B, 8.5 and 8.4 Corollary}, it is
a diagonalizable group.  Let $K'$ be a finite extension of $K$ over which
$C_K$ splits, and let $S'$ be a corresponding generically finite cover
of $S$.  Then, by \cite{B, 8.7}, $C_K\times_K K'$ is of the form
$F\times\Gm^\mu$, where $F$ is a finite group and $\mu\in\Bbb N$.
Hence there exists a diagonalizable group $C$ over $k$ and a nonempty
open subgroup $U$ of $S'$ such that the closure of $B_K\times_K K'$
in $\Cal B\times_S U$ is $U$\snug-isomorphic to $C\times_k U$.
By rigidity \cite{B, 8.10}, it follows that the induced map $C\times_k U\to A$
factors through the projection onto the first factor; hence we may regard $C$
as a subgroup of $A$.  Since $\Cal B$ is closed, it follows that
$\Cal B\supseteq C\times_k S$.

After replacing $A$ with $A/C$ and $\Cal B$ with its image in $(A/C)\times_k S$,
we may assume that $\Cal B$ is generically finite over $A_0\times_k S$,
of degree $1$.  Since $\Cal B$ is also reduced, it corresponds to a (reduced)
rational point on the generic fiber of $\rho_{(S)}\restrictedto{\Cal B}$,
hence $\Cal B$ is the closure of the image of a rational section
$\sigma\:U\to A\times_k S$,
where $U$ is an open dense subset of $A_0\times_k S$.  Translating by closed
points of $A_0$ and using the fact that $\Cal B$ is a group subscheme,
we see that $U$ is of the form $A_0\times_k V$ for some open dense
subset $V\subseteq S$ (and that $\sigma$ is a homomorphism of $V$\snug-group
schemes).  Thus $\Cal B$
is a family of regular sections of $\rho\:A\to A_0$, parametrized by $V$.
But $\rho\:A\to A_0$ admits at most one regular section passing through
the group identity of $A$, since the ratio of any two such sections is a
regular map $A_0\to\Gm^{\mu(A)}$.  Thus $\Cal B\cap (A\times_k V)=B'\times_k V$
for some $B'\subseteq A$.  Since $\Cal B$ is closed, it follows
that $\Cal B\supseteq B'\times_k S$.  Since all irreducible components of
$\Cal B$ dominate $A_0\times_k S$, we have $\Cal B=B'\times_k S$.
Finally, $B'=\Cal B_s=B$, and we are done.\qed
\enddemo

Recall that a {\bc regular} field extension is a field extension $K/k$
such that $K$ is linearly disjoint from $\bar k$ over $k$; equivalently
(in characteristic zero), $k$ is algebraically closed in $K$.

\lemma{\05.14}  Let $k$ be any field, let $K/k$ be a regular field extension,
let $A$ be a geometrically integral semiabelian variety over $k$, and
let $B$ be a geometrically integral group subvariety of $A_K:=A\times_k K$.
Then there exists a group subvariety $B_k$ of $A$ such that $B=B_k\times_A A_K$.
\endit

\demo{Proof}  This lemma is already known when $A$ is an abelian variety:
see \cite{Mi, Cor.~20.4}.  The proof here is essentially the same proof,
using the stronger Lemma \05.13.

Let $S$ be a variety over $k$ with $K(S)=K$.  Since $K/k$ is regular,
$S$ is geometrically integral.  Let $\Cal B$ be the closure of $B$
in $A\times_k S$.  After replacing $S$ with an open subvariety, 
we may assume that $\Cal B$ is smooth over $S$,
that $\rho_{(S)}\restrictedto{\Cal B}\:\Cal B\to A_0\times_k S$ is smooth,
and that the morphisms $\Spec K\to B$, $B\to B$, and $B\times_K B\to B$
expressing the group structure of $B$ extend to morphisms $S\to\Cal B$, etc.
This latter condition ensures that $\Cal B$ is a group subscheme of
$A\times_k S$.

By \cite{R, p\. 412}, there is a finite separable extension $k'$ of $k$
such that $S(k')\ne\emptyset$.  We may assume that $k'$ is Galois over $k$.
After base change to $k'$, Lemma \05.13 implies the existence of
$B_{k'}\subseteq A\times_k k'$ such that $B\times_K Kk'=B_{k'}\times_{k'}Kk'$
(as subvarieties of $A\times_k Kk'$).  This equality determines $B_{k'}$
uniquely; hence it is invariant under $\operatorname{Gal}(k'/k)$.
By descent (see, for example, \cite{Se~1, V.20}), $B_{k'}$ is of the form
$B_k\times_k k'$ for some semiabelian subvariety $B_k$ of $A$.\qed
\enddemo

\lemma{\05.15}  Let $X$ be a closed subvariety of a semiabelian variety $A$.
Then $\bar\kappa(X)\ge0$\snug, with equality if and only if $X$ is a
translated semiabelian subvariety of $A$.
\endit

\demo{Proof}  See \cite{I~2, Thm.~4 and Thm.~2}.\qed
\enddemo

\thm{\05.16}  Let $A$ be a semiabelian variety, let $X$ be a closed subvariety
of $A$, and let $D$ be an effective Weil divisor on $X$.  Then
\roster
\myitem a.  $\bar\kappa(X\setminus D)\ge 0$, with equality if and only if
$X$ is a translated semiabelian subvariety of $A$ and $D=0$;
\myitem b.  $\bar\kappa(X\setminus D)+\dim B(X\setminus D)=\dim X$; and
\myitem c.  $B(X\setminus D)=0$ if and only if $X\setminus D$ is of
logarithmic general type.
\endroster
\endit

\demo{Proof}  The inequality $\bar\kappa(X\setminus D)\ge0$ is immediate
from Lemma \05.15 and the inequality $\bar\kappa(X\setminus D)\ge\bar\kappa(X)$.
The remainder of part (a) follows by Lemmas \05.15 and \05.12.

Next consider (b).  By Hironaka's resolution of singularities and by
\cite{I~3, Thm. 10.3}, there exists a complete nonsingular variety $X'$,
a normal crossings divisor $D'$ on $X'$, a proper birational morphism
$\pi\:X'\setminus D'\to X\setminus D$, a complete nonsingular variety $W$,
and a dominant morphism $\Phi\:X'\to W$ such that
\roster
\myitem i.  $\dim W=\bar\kappa(X\setminus D)$,
\myitem ii.  $K(X')$ is regular over $K(W)$, and
\myitem iii.  if $X'_\eta$ denotes the generic fiber of $\pi$,
then $h^0(X'_\eta,mK_{(X',D')})\le1$ for all $m\in\Bbb N$, and is not
always zero.
\endroster
The canonical divisor of the generic fiber of $\pi$ is
$K_{X'}\restrictedto{X'_\eta}$; this follows by looking at the First
Exact Sequence for differentials, which is exact on the left in this case,
and taking highest exterior powers.  Therefore the logarithmic canonical
divisor of the generic fiber of $\pi\restrictedto{X'\setminus D'}$
is $(K_{X'}+D')\restrictedto{X'_\eta}$.  Thus, by (iii), the logarithmic
Kodaira dimension of the generic fiber of $\pi\restrictedto{(X'\setminus D')}$
is zero.  We may regard this as a subvariety of $A_\eta:=A\times_k K(W)$;
after base change to $\widebar{K(W)}$, it is a translated semiabelian
subvariety of $A\times_k\widebar{K(W)}$.  By Lemma \05.14, that subvariety
comes from a subvariety $B$ of $A$.  Since the generic fiber is invariant
under translation by $B$, the same holds for $X\setminus D$.  Thus
$B(X\setminus D)\supseteq B$, and therefore by (i),
$$\dim B(X\setminus D) \ge \dim B = \dim X - \bar\kappa(X\setminus D)\;.
  \tag\05.16.1$$
But also, by Lemmas \05.6 and \05.7,
$$\bar\kappa(X\setminus D)
  = \bar\kappa\bigl((X\setminus D)/B(X\setminus D)\bigr)
  \le \dim X - \dim B(X\setminus D)\;.$$
Combining this with (\05.16.1) gives part (b).

Part (c) follows immediately from (b).\qed
\enddemo

\beginsection{\06}{The Kawamata Structure Theorem}

This section generalizes the Kawamata Structure Theorem to the context
of a closed subvariety of a semiabelian variety, minus a divisor.

\thm{\06.1}  Let $A$ be a semiabelian variety over an algebraically closed
field of characteristic zero, let $X$ be a closed subvariety of $A$,
and let $D$ be an effective divisor on $X$.
Let $Z=Z(X\setminus D)$ be the union of all positive dimensional
translated semiabelian subvarieties of $A$ contained in $X\setminus D$.
Then $Z$ is a Zariski-closed subset of $X\setminus D$, and each
irreducible component has nontrivial Ueno fibration.
\endit

\demo{Proof}
We may assume that $X\setminus D$ has trivial Ueno fibration,
for otherwise the theorem is trivial (with $Z=X\setminus D$).
By noetherian induction it then suffices to show that $Z$ is not Zariski-dense.
Let $B=B(X)$ and $X'=X/B$; then $X'$ has trivial Ueno fibration and
there is a fiber bundle $\theta\:X\to X'$ with fiber $B$.

Let $Z'$ (resp\. $Z''$) be the union of all nontrivial translated
semiabelian subvarieties of $A$ which are contained in $X\setminus D$
and which lie (resp\. do not lie) in fibers of $\theta$.  Then $Z=Z'\cup Z''$.
But $Z''\subseteq\theta^{-1}(Z(X'))$, which is not Zariski-dense by the
Kawamata Structure Theorem for {\it closed\/} subvarieties of
semiabelian varieties \cite{N, Lemma~4.1}.
Thus it suffices to show that $Z'$ is not Zariski-dense.

Let $Z_0$ be the union of all nontrivial translated {\it abelian\/}
subvarieties making up $Z'$, and let $Z_1$ be the union of all nontrivial
images of $\Gm$ lying in $X\setminus D$ and contained in fibers of $\theta$.
Then
$$Z'=Z_0\cup Z_1\;.$$
(Note that a nontrivial image of $\Gm$ is a translated semiabelian subvariety,
by \cite{I~2, Thm.~2}.)

First consider $Z_0$.  Suppose $Z_0\ne\emptyset$ and let $C$ be a translated
abelian subvariety of $A$ lying in $X\setminus D$ and lying in a fiber
of $\theta$.  Then the closure of $D$ in $\widebar A$ does not meet $C$;
hence this is true of all translates of $C$ unless they lie in the closure
of $D$.  Thus $B(X\setminus D)$ contains the abelian subvariety corresponding
to $C$.  This contradicts the assumption that $B(X\setminus D)$ is trivial,
so $Z_0=\emptyset$.

This leaves $Z_1$.
The remainder of this proof is motivated by the proof of Theorem \00.2.

Let $B_0$ be the abelian quotient of $B$ and let $B_1$ be the toric part,
so that there exists an exact sequence $0\to B_1\to B\to B_0\to0$.
Then $\theta\:X\to X'$ factors as
$$ X @>\phi>> X_0 @>\psi>> X'\; ,$$
where $\phi$ and $\psi$ are fiber bundles with fibers $B_1$ and $B_0$,
respectively.

Let $\pi_0\:\widetilde X_0\to X_0$, $X\subseteq\widebar X$,
$\pi\:\widetilde X\to\widebar X$, and $\widetilde D$ be as in Theorem \02.4.
Let $\widebar B_1$ be the equivariant completion of $B_1$
corresponding to $\widebar X$.  Let $\chi\:\Gm\to X\setminus D$ be a nontrivial
morphism whose image is contained in a fiber of $\theta$.  Then its image must
lie in a fiber of $\phi$, say, $\operatorname{Im}\chi\subseteq \phi^{-1}(x)$,
$x\in X_0$.  Let $\tilde x\in\widetilde X_0$ be a point lying over $x\in X_0$;
then $\chi$ lifts to $\widetilde\chi\:\Gm\to\widetilde X$.
Since $\widetilde D$ is ample on fibers
of $\tilde \phi\:\widetilde X\to\widetilde X_0$, the closure of the image
of $\widetilde\chi$ must meet $\widetilde D$.  Let $P$ be a point where
they meet, let $T$ be the $B_1$\snug-orbit of $\widebar B_1$ corresponding
to the $B_1$\snug-orbit of $\tilde \phi^{-1}(\tilde x)$ containing $P$,
and let $p\:U\to T$ be the projection defined in Proposition \02.5;
$U\supseteq B_1$.  This $p$ defines a projection $q$ from an open subset
of $\widetilde X$ to the subset $\widetilde T$ of $\widetilde X$
corresponding to $T$.  Then $q\circ\widetilde\chi$ defines a morphism
$\Bbb A^1\to T$.  Since $T$ is isomorphic to a product of $\Gm$'s,
it follows that this morphism must be trivial.
Thus the image of $\widetilde\chi$ must lie in the proper Zariski-closed
subset $q^{-1}(\Supp\widetilde D\cap\widetilde T)$, so the image
of $\chi$ must lie in the corresponding proper Zariski-closed subset of $X$.

Since there are only finitely many such $\widetilde T$,
it follows that $Z_1$ cannot be Zariski-dense.\qed
\enddemo

\enddocument